\documentclass[12pt,twoside]{amsart}
\usepackage{amsfonts,amsmath}
\usepackage{pinlabel}
\usepackage{enumerate}
\usepackage{mathtools}
\usepackage{xcolor,colortbl}
\usepackage{url}

\def\endpf{\relax\ifmmode\expandafter\endproofmath\else
  \unskip\nobreak\hfil\penalty50\hskip.75em\hbox{}\nobreak\hfil\bull
  {\parfillskip=0pt \finalhyphendemerits=0 \bigbreak}\fi}
\def\bull{\vbox{\hrule\hbox{\vrule\kern3pt\vbox{\kern6pt}\kern3pt\vrule}\hrule}}

\newtheorem{defn}{Definition}[section]
\newtheorem{lemma}[defn]{Lemma}

\newtheorem{theorem}[defn]{Theorem}

\newtheorem{proposition}[defn]{Proposition}

\newtheorem{maintheorem}{Theorem}
\newtheorem{mainproposition}[maintheorem]{Proposition}
\newtheorem{maincorollary}[maintheorem]{Corollary}
\newtheorem{example}[defn]{Example}

\newcommand{\zz}{{\mathbb Z}}

\newcommand{\qq}{{\mathbb Q}}

\newcommand{\ozsvath}{Ozsv\'{a}th}
\newcommand{\szabo}{Szab\'{o}}
\newcommand{\spin}{\ifmmode{\rm Spin}\else{${\rm spin}$\ }\fi}
\newcommand{\spinc}{\ifmmode{{\rm Spin}^c}\else{${\rm spin}^c$\ }\fi}
\newcommand{\spinct}{\mathfrak t}
\newcommand{\spincs}{\mathfrak s}
\newcommand{\tors}{{\it Tors}}

\newcommand{\calt}{\mathcal{T}}

\newcommand{\tdelta}{\widetilde\Delta}
\newcommand{\tw}{\widetilde W}
\newcommand{\tf}{\widetilde F}
\newcommand{\td}{\tilde d}

\newcommand{\CP}{\mathbb{CP}}

\definecolor{Gray}{gray}{0.8}
\newcommand{\gc}{\cellcolor{Gray}}

\newenvironment{narrow}[2]{%
 \begin{list}{}{%
  \setlength{\topsep}{0pt}%
  \setlength{\leftmargin}{#1}%
  \setlength{\rightmargin}{#2}%
  \setlength{\listparindent}{\parindent}%
  \setlength{\itemindent}{\parindent}%
  \setlength{\parsep}{\parskip}%
 }%
\item[]}{\end{list}}

\newif\ifpic
\picfalse    
\pictrue   

\DeclareMathOperator{\Char}{Char}

\begin{document}

\title{Immersed disks, slicing numbers and concordance unknotting numbers }
\author[Brendan Owens]{Brendan Owens}
\address{School of Mathematics and Statistics \newline\indent 
University of Glasgow \newline\indent 
Glasgow, G12 8QW, United Kingdom}
\email{brendan.owens@glasgow.ac.uk}
\author[Sa\v so Strle]{Sa\v so Strle}
\address{Faculty of Mathematics and Physics \newline\indent 
University of Ljubljana \newline\indent Jadranska 21 \newline\indent 
1000 Ljubljana, Slovenia }
\email{saso.strle@fmf.uni-lj.si}
\date{\today}
\thanks{B. Owens was supported in part by  EPSRC grant EP/I033754/1.}
\thanks{S. Strle was supported in part by ARRS Grant P1--0288}

\begin{abstract}  We study three knot invariants related to smoothly immersed disks in the four-ball.  These are the four-ball crossing number, which is the minimal number of normal double points of such a disk bounded by a given knot; the slicing number, which is the minimal number of crossing changes required to obtain a slice knot; and the concordance unknotting number, which is the minimal unknotting number in a smooth concordance class.
Using Heegaard Floer homology we obtain bounds that can be used to determine two of these invariants for all prime knots with crossing number ten or less, and to determine the concordance unknotting number for all but thirteen of these knots.  As a further application we obtain some new bounds on Gordian distance between torus knots.
We also give a strengthened version of \ozsvath\ and \szabo's obstruction to unknotting number one.
\end{abstract}

\maketitle

\pagestyle{myheadings}
\markboth{BRENDAN OWENS AND SA\v{S}O STRLE}{}


\section{Introduction}
\label{sec:intro}

The unknotting number $u(K)$ of a knot $K$ in $S^3$ is the minimal number of crossing changes required to convert it to an unknot. The trace of a regular homotopy realizing the crossing changes yields a normally immersed annulus $A$ in $S^3\times [0,1]$ with a singularity for each crossing change. We say a surface is \emph{normally immersed} if the immersion is proper (sends precisely boundary to boundary) and the only singularities are normal double points (also known as normal crossings), that is to say, transverse double points in the interior of the surface. Since the other boundary of $A$ is an unknot, we can complete $A$ to a normally immersed disk $\Delta$ in $B^4$ with boundary $K$ and $u(K)$ double points. Minimising the number of double points in any such disk with boundary $K$ gives a concordance invariant of $K$, the \emph{$4$-ball crossing number} $c^*(K)$   \cite{tk,my,shibuya}, also referred to as the $4$-dimensional clasp number. Recall that $K$ and $K'$ are \emph{concordant} if they cobound a properly embedded annulus in $S^3 \times I$. A knot $K$ is \emph{slice} if it is concordant to the unknot or equivalently if it bounds a smoothly embedded disk in $B^4$ -- such a disk is called a \emph{slice disk} or a \emph{nullconcordance}.

We show in Proposition \ref{prop:movie} that any normally immersed disk in $B^4$ can be factored into a concordance, followed by the trace of a regular homotopy, followed by a nullconcordance; thus it is natural to consider two intermediate invariants. The first is the \emph{slicing number} $u_s(K)$ \cite{ask,liv,slicing,rudolph,shibuya}, which is the minimal number of crossing changes required to obtain a slice knot. A set of such crossing changes clearly gives rise to a normally immersed disk, obtained by capping off the trace of the crossing-change homotopy by a slice disk. 

The second intermediate invariant is the \emph{concordance unknotting number} $u_c(K)$, which is the minimal unknotting number of any knot in the concordance class of $K$. 
Noting that a normal double point may be resolved at the cost of increasing the genus of the immersed surface by one, we have the following inequalities:
\begin{equation}
\label{eq:ineq}
u(K) \ge u_c(K), u_s(K) \ge c^*(K) \ge g^*(K) \ge |\sigma(K)|/2,
\end{equation}
where $g^*(K)$ denotes the smooth four-ball genus of $K$, $\sigma(K)$ denotes the signature, and the last inequality is due to Murasugi \cite{murasugi}.   In this article we develop some tools to calculate these invariants, and we determine each of $u_s$ and $c^*$ for all prime knots with 10 or fewer crossings, and $u_c$ for all but three knots up to 9 crossings and all but ten 10-crossing knots.  Our results make use of Montesinos' trick \cite{mont}, which implies that the double branched cover of an unknotting number one knot is given by a half-integer surgery on a knot in $S^3$, and also theorems of 
Cochran and Lickorish \cite{cl-unknotting},  \ozsvath\ and \szabo\ \cite{osu1}, and Ni and Wu \cite{niwu}.

Theorem \ref{thm:imm} is a refinement of \cite[Theorem 3.7]{cl-unknotting} concerning a four-manifold bounded by the double branched cover of a knot and determined by an immersed disk bounded by the knot. 
From it we obtain the following result, which shows that the slicing number obstruction given in \cite[Theorem 2]{slicing} in fact applies to the four-ball crossing number.

\begin{maintheorem}
\label{thm:clasp}
Suppose that $K$ bounds a normally immersed disk with $r_+$ positive double points and $r_-=\sigma(K)/2$ negative double points.
Then the branched double cover $\Sigma(K)$ bounds a positive-definite smooth four-manifold $X$ with $b_2(X) = 2(r_++r_-)$ whose intersection form $Q_X$ is of half-integer surgery type, with exactly $r_+$ of the diagonal entries odd, and $\det Q_X$ divides $\det K$ with quotient a square.
\end{maintheorem}

Recall that a quadratic form over the integers is of \emph{half-integer surgery type} if with respect to some basis for the lattice it is represented by a matrix of the form
$$\left[\begin{matrix} A & I \\ I & 2I\end{matrix}\right],$$
where $I$ denotes the identity matrix.

In case of $c^*(K)=1$, we obtain an explicit obstruction
in terms of the correction terms of the double branched cover of $K$. 
Recall that the set of \spinc structures on a three-manifold $Y$ is an affine copy of $H_1(Y;\zz)$, with additional structure given by conjugation of \spinc structures, fixed points of which correspond to spin structures on $Y$.  When $Y$ is the double cover of $S^3$ branched along a knot, there is a canonical identification of $\spinc(Y)$ with the homology group of $Y$ by using the unique spin structure as the origin; we use this implicitly throughout.  Recall also that for a \spinc structure $\spincs$ on a rational homology three-sphere $Y$, \ozsvath\ and \szabo\ \cite{os4} defined the correction term or $d$-invariant, $d(Y,\spincs)$, as the absolute grading of a particular subgroup of the Heegaard Floer homology group of $(Y,\spincs)$.  The $d$-invariants are rational numbers which are computable in many cases.

\begin{maincorollary}\label{cor:c=1}
Let $K\subset S^3$ be a knot with four-ball crossing number one. Suppose $K$ has signature two and let $Y=\Sigma(K)$ be the double branched cover of $K$. Then for some factorization $\det K=rs^2$ there exists an order $rs$ subgroup $H \le \spinc(Y)$ and an epimorphism $\phi \colon H \to \zz/r\zz$ such that the normalized $d$-invariants
$$\td_t=-d(Y,t)+i^2/2r-\begin{cases}0& i \equiv (r-1)/2 \pmod 2\\1/2&i \equiv (r+1)/2 \pmod 2\end{cases}$$
for $i=0,\dots,(r-1)/2$ and for every $t\in \phi^{-1}(i)$ satisfy
\begin{enumerate}[(i)]
\item \label{itemc1} positivity: $\td_t\ge0$;
\item \label{itemc2} evenness: $\td_t\in 2\, \zz$.
\end{enumerate}
If there exists an epimorphism $\phi$ as above, we say that $Y$ admits a \emph{positive even subgroup matching}.

If $K$ has signature zero, then at least one of $\pm\Sigma(K)$ admits a positive even subgroup matching.
\end{maincorollary}

Note that by conjugation invariance of $d$-invariants, the constraint given in Corollary 2 extends to $-(r-1)/2\le i \le (r-1)/2$.  This applies to Theorems \ref{thm:u=1} and \ref{thm:uc=1} below as well.

The following theorem concerning manifolds given as half-integer surgery on a knot in $S^3$ is an extension of results of \ozsvath\ and \szabo\ \cite{osu1}, where the first three conditions are established.

\begin{maintheorem}
\label{thm:u=1}
Suppose that a rational homology sphere $Y$, with $|H_1(Y)|=r$ odd, is given by $r/2$ surgery on a knot  in $S^3$.  Then there exists a group isomorphism $$\phi:\spinc(Y)\to\zz/r\zz$$ such that the normalized $d$-invariants
$$\td_i=-d(Y,\phi^{-1}(i))+i^2/2r-\begin{cases}0& i \equiv (r-1)/2 \pmod 2\\1/2&i \equiv (r+1)/2 \pmod 2\end{cases}$$
for $i=0,\dots,(r-1)/2$ satisfy the following conditions:
\begin{enumerate}[(i)]
\item \label{item1} positivity: $\td_i\ge0$;
\item \label{item2} evenness: $\td_i\in 2\, \zz$;
\item \label{item3} symmetry: 
$$\td_{2j}=\td_{2j+1}\ \text{for}\ 0\le j\le (r-5)/4\quad \text{if}\ r\equiv 1 \pmod 4\text{,}$$
$$\td_{2j-1}=\td_{2j}\ \text{for}\ 1\le j\le (r-3)/4\quad \text{if}\ r\equiv -1 \pmod 4\text{;}$$
\item \label{item4} monotonicity: $\td_{i}\le\td_{i+1}$ for $0\le i<(r-1)/2$;
\item \label{item5} boundedness: $\td_{i+1}\le\td_i + 2$ for $0\le i<(r-1)/2$.
\end{enumerate}
\end{maintheorem}

If there exists an isomorphism $\phi$ satisfying the conditions of Theorem \ref{thm:u=1} we say that $Y$ admits a \emph{positive even symmetric monotone matching}.  
In case the manifold $Y$ is an $L$-space the normalized $d$-invariants $\td_i$ of Theorem \ref{thm:u=1} are given by the torsion coefficients defined using the Alexander polynomial of the surgery knot, as in Theorem 1.2 of \cite{Qsurg}.  Conditions \eqref{item4} and \eqref{item5} follow easily in that case.

For many examples previously obstructed using the symmetry condition, the monotonicity condition may be substituted.  In fact, we have not yet found an example which is obstructed by symmetry but not by monotonicity, or vice versa.
More interestingly perhaps, the obstruction given by conditions \eqref{item1}, \eqref{item2} and \eqref{item4} of Theorem \ref{thm:u=1} extends over rational homology cobordisms as in the following statement. 

\begin{maintheorem}
\label{thm:uc=1}
Suppose that a rational homology sphere $Y$, with $|H_1(Y)|$ odd, is rational homology cobordant to positive half-integer surgery on a knot in $S^3$.  Then for some factorisation $|H_1(Y)|=rs^2$ there exists an order $rs$ subgroup $H \le \spinc(Y)$ and an epimorphism $\phi \colon H \to \zz/r\zz$ such that the $d$-invariants of $Y$ are constant on the fibers of $\phi$.  
Moreover, the normalized $d$-invariants
$$\td_i=-d(Y,\phi^{-1}(i))+i^2/2r-\begin{cases}0& i \equiv (r-1)/2 \pmod 2\\1/2&i \equiv (r+1)/2 \pmod 2\end{cases}$$
for $i=0,\dots,(r-1)/2$ satisfy the following conditions:
\begin{enumerate}[(i)]
\item \label{uc1}positivity: $\td_i\ge0$;
\item \label{uc2}evenness: $\td_i\in 2\zz$;
\item \label{uc3}monotonicity: $\td_{i}\le\td_{i+1}$ for $0\le i<(r-1)/2$.
\end{enumerate}
\end{maintheorem}

Here we are using $d(Y,\phi^{-1}(i))$ to denote the constant value of the $d$-invariant on the fiber $\phi^{-1}(i)$ for $i\in\zz/r\zz$.
If there exists an epimorphism $\phi$ satisfying the conditions  
of Theorem \ref{thm:uc=1} we say that $Y$ admits a \emph{positive even monotone subgroup matching}.

Applying these results to knots and keeping track of signs we obtain

\begin{maincorollary}\label{cor:knots}
Let $K\subset S^3$ be a knot with unknotting number one.  If $K$ has signature two, then the double branched cover $\Sigma(K)$ admits a positive even symmetric monotone matching.  If $K$ has signature zero, then at least one of $\pm\Sigma(K)$ admits a positive even symmetric monotone matching.

Suppose $K$ has concordance unknotting number one. If $K$ has signature two, then $\Sigma(K)$ admits a positive even monotone subgroup matching.
If $K$ has signature zero, then at least one of $\pm\Sigma(K)$ admits a positive even monotone subgroup matching.
\end{maincorollary}

For each of the inequalities in \eqref{eq:ineq} except for $u_s\ge c^*$ there exist examples for which the inequality is strict.  One may ask whether in fact $u_s(K)=c^*(K)$ for all knots.  We will see in Section \ref{sec:examples} that this equality holds for all prime knots with 10 or fewer crossings.
Another question which seems to be open is whether the slicing number $u_s$ is a concordance invariant.  
It turns out these questions are related to a generalisation of Fox's  slice-ribbon question, which asks whether every slice knot is in fact ribbon (admits a slice disk which is ribbon).  
A properly embedded or immersed surface in $B^4$ is called \emph{ribbon} if the restriction of the radial distance function is Morse without local maxima.

\begin{mainproposition}\label{prop:usconc}
Let $c^*_r(K)$ denote the minimal number of double points in a normally immersed ribbon disk bounded by $K$ in $B^4$, and let {\rm SRC} be the Slice-Ribbon Conjecture, which states that all slice knots are ribbon.  Then:
\begin{align*}
c^*=c^*_r&\iff \mathrm{SRC} \text{ and }u_s=c^*\\
\intertext{and moreover}
u_s=c^*&\iff u_s \text{ is a concordance invariant}.
\end{align*}
\end{mainproposition}

In the last two sections of the paper we consider examples.  We compute the slicing number and four-ball crossing number for all prime knots of ten crossings or fewer, and the concordance unknotting number for all but thirteen of these knots.  We also obtain some new bounds on Gordian distance between torus knots.

\vskip2mm
\noindent{\bf Acknowledgements.}  
We are grateful to Josh Greene for a helpful comment about $L$-spaces which led us to condition \eqref{item5} in Theorem \ref{thm:u=1}, and to  Frank Swenton who helped us to use his Kirby Calculator software \cite{klo} to produce and manipulate knot diagrams.  We thank Maciej Borodzik who suggested we look at Gordian distances between torus knots.  We thank the anonymous referees for helpful suggestions to improve the exposition.


\section{Geometric constructions}
\label{sec:geometric}

In this section we collect some results regarding normally immersed surfaces in $B^4$, and prove Theorem \ref{thm:clasp}, Corollary \ref{cor:c=1}, and Proposition \ref{prop:usconc}.

\begin{figure}[htbp]
\begin{center}
\ifpic{
\labellist
\small\hair 3pt
\pinlabel $\longrightarrow$ at 120 80
\endlabellist
\centering
\includegraphics[scale=0.6]{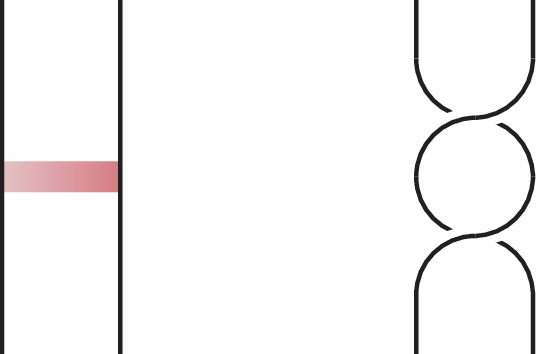}
}
\else \vskip 5cm \fi
\begin{narrow}{0.3in}{0.3in}
\caption
{\bf Encoding a crossing change with a band.}
\label{fig:redband}
\end{narrow}
\end{center}
\end{figure}

Recall that a crossing change in a link $L$ may be recorded by placing a framed arc or equivalently a band connecting two arcs of the link.  The
crossing change operation consists of replacing the two arcs of the band on $L$ with a full positive twist as in Figure \ref{fig:redband}.

\begin{proposition}
\label{prop:movie}
Let $F$ be a  connected surface normally immersed in $S^3 \times I$ with $L=F \cap (S^3 \times \{0\})\ne \emptyset$. After an isotopy rel boundary we may assume that $F \cap (S^3 \times [0,1/3])$ is a concordance, $F \cap (S^3 \times [1/3,2/3])$ is the trace of a regular homotopy, and $F \cap (S^3 \times [2/3,1])$ is a smoothly embedded ribbon surface (that is, the projection to the interval on this part of the surface is a Morse function without local minima). 

In particular, if $F$ is an immersed disk, then it can be factored into a concordance, followed by the trace of a regular homotopy, followed by a ribbon nullconcordance.
\end{proposition}

\proof
After a preliminary isotopy we may assume the projection to the interval is a Morse function $h$ on $F$ whose critical points are distinct from the double points.  Then generic level sets $F\cap (S^3 \times \{t\})$ are smooth links in the three-sphere, and a sequence of diagrams of such links gives a ``movie presentation'' of the surface $F$.
Neighborhoods of double points of $F$ in the movie presentation of $F$ relative to $h$ correspond to crossing changes. Successive frames in the movie picture of $F$ are then obtained by crossing changes, Morse moves and isotopies. Finally, after an isotopy of the surface supported in a small neighborhood of all index zero and two critical points, we may arrange that all the minima of $h$ occur before and all the maxima occur after all the saddle and double points. We let $F' \subset F$ be a subsurface bounded by two regular level sets that contains all the saddle and double points of $F$, but no minima and maxima.

A complete description of the surface $F'$ can be encoded in a single diagram of the link $L'$ which is the lower boundary of $F'$; thus $L'$ consists of $L$ and the unlink of co-attaching circles for handles of index zero. Following the movie we encode each crossing change encountered by adding a red band as in Figure \ref{fig:redband} and each addition of a handle of index one (a band move) by adding a blue one. By shortening the arcs and pulling the rest of $L'$ along we can construct a diagram for $L'$ relative to which the added colored bands lie in a plane, are nonoverlapping, and intersect $L'$ only in attaching arcs. It is clear from this picture that we can now implement the crossing changes and band moves in any order we desire.

Since $F$ is connected we can first modify $L'$ by one blue band move (1-handle addition) for each of the components corresponding to the boundaries of handles of index zero, so that each becomes connected to a component of $L$.   This forms a concordance, which is the first stage of factorisation. The second stage of factorisation consists of all the crossing changes encoded by the red bands. The remaining handles of index one (blue bands) and two form a ribbon surface which is the last stage of factorisation.
\endproof

It follows from the above proposition that the four-ball crossing number of a knot is equal to the minimal slicing number of any knot in its concordance class:

\begin{equation}\label{eq:concus}
c^*(K)=\min_{K'\text{ concordant to } K}u_s(K').
\end{equation}

\proof[Proof of Proposition \ref{prop:usconc}]
A knot $K$ is slice if and only if $c^*(K)=0$, and is ribbon if and only if  $c^*_r(K)=0$; thus $c^*=c^*_r$ implies SRC.  It was observed by Shibuya \cite{shibuya} and Rudolph \cite{rudolph} (and follows easily from 
Proposition \ref{prop:movie}) that $c^*_r(K)$ is equal to $u_r(K)$ which is the minimal number of crossing changes from $K$ to a ribbon knot. Note that SRC implies $u_s=u_r$.  The equivalence
$$c^*=c^*_r\iff  \mathrm{SRC} \text{ and }u_s=c^*$$
now follows by considering equalities among the quantities $c^*$, $c^*_r$, $u_r$ and $u_s$.
Finally, the nontrivial part of the equivalence
$$u_s=c^*\iff \text{ $u_s$ is a concordance invariant}$$
follows easily from Proposition \ref{prop:movie} via \eqref{eq:concus}.
\endproof

The following theorem implies Theorem \ref{thm:clasp} and may be used to give a new proof of  \cite[Theorem 2]{slicing}, noting that the trace of a homotopy given by a crossing change sequence from $K$ to a slice knot $J$ may be glued to a slice disk for $J$ to give an immersed disk $\Delta$ in $B^4$ bounded by $K$.  There is a sign error in Remark 3.5 of that paper: changing a positive (respectively negative) crossing in $K$ results in a positive (respectively negative) double point of $\Delta$. 

\begin{theorem}
\label{thm:imm}
Let $\Delta$ be a normally immersed disk in the four-ball with $r$ double points and boundary $K$. Then the double branched cover $\Sigma(K)$ of $K$ bounds a manifold with $b_2=2r$ and signature $\sigma(K)+2r_+$, where $\sigma(K)$ denotes the signature of $K$ and $r_+$ the number of positive double points in $\Delta$. Moreover, the intersection pairing of this manifold is of half-integer surgery type and the number of odd squares in such a representation of the pairing is equal to $r_+$. 
\end{theorem}

\proof 
We follow the construction in \cite[Theorem 3.7]{cl-unknotting} where all the properties except those in the last sentence of the theorem are established. We briefly recall the construction. Blow up $B^4$ at every double point of $\Delta$ to obtain $W=B^4 \# r\CP^2$ with an embedded disk $\tdelta$ and let $\tw$ be the double branched cover of $W$ with branch set $\tdelta$. Then $b_1(\tw)=0$ and hence $b_2(\tw)=2r$ via an Euler characteristic computation, and the signature formula follows from the $G$-signature theorem, as in the proof of \cite[Theorem 3.7]{cl-unknotting}. 
We exhibit below a collection of $2r$ independent homology classes in $H_2(\tw;\zz)$ with the claimed type of pairing. 
It follows then from \cite[Proposition 2.4]{slicing} that the pairing on $\tw$ is of the same type.

Let $p$ be a double point of $\Delta$ and let $B$ be a small ball around $p$; we may assume $\Delta$ intersects the boundary of $B$ in two great circles. Denote by $E$ the total space of the degree one disk-bundle over $S^2$ that replaces $B$ after the blow-up at $p$. The proper transform of $\Delta$ thus contains two fibres of $E$ and in the double branched cover $\tw$ the zero section of $E$ lifts to a sphere $S$ of self-intersection two. 

Choose a loop $C$ in $\Delta$ that is the image of an arc connecting the two points in the disk mapping to $p$ under the immersion to $\Delta$. We may assume $C$ does not contain any other double point of $\Delta$. Choose a (short) vector $v_p$ transverse to both sheets of $\Delta$ at $p$ and extend it to a (short enough) nonvanishing normal vector field $v$ to $\Delta$ along $C$. Then $C$ and its pushoff $C_v$ along $v$ cobound an embedded annulus $A$. Since the first homology of the complement of $\Delta$ is generated by the meridian $\mu$ of $\Delta$, it follows that $C_v$ is homologous to $k[\mu]$ for some $k\in \zz$. Choose $|k|$ disjoint curves $L\subset \Delta - B$, parallel to a component of $\Delta \cap \partial B$, and let $F_0$ be an embedded surface bounded by $C_v \cup L$ that intersects $\Delta$ only in $L$. Then $F_0 \cup A$ is a surface (with possible self-intersections) with boundary $C \cup L$. We may replace it by a surface $F$ with the same boundary and only (transverse) self-intersections along the boundary. 

In the blow-up process we replace $B$ by $E$ and $F$ by $F'=F - \operatorname{int} B$. By making $B$ smaller if necessary we may further assume that $F' \cap \partial E$ is a section $s_0$ of $E$ over an arc in the base sphere connecting the branch points. Extend $s_0$ to a section $\Sigma$ of $E$ that intersects the zero section transversely in one point. Then the preimage $\tf'$ of $F'$, in the branched double cover $\tw$, is a normally immersed surface with boundary the preimage $\widetilde{s_0}$ of $s_0$, since it is the double of $F'$ along the part of its boundary that lies in $\Delta$. Similarly, $\widetilde{s_0}$ separates the preimage $\widetilde \Sigma$ of $\Sigma$ into two ``hemispheres''. Hence $\tf'$ and one hemisphere of $\widetilde \Sigma$ form a closed normally immersed surface $\tf$ in $\tw$. By construction the intersection number of $\tf$ with $S$ is $\pm 1$, so they represent a dual pair of homology classes. Since the sphere $S$ is contained in the preimage of $E$, spheres corresponding to different double points clearly represent different classes which are disjoint.  Thus we obtain $2r$ homology classes with the required form of intersection pairing.

To determine the parity of the self-intersection of $\tf$ construct a pushoff as follows. Choose a short nonvanishing normal vector field $u$ along each component of $\partial F$ in $\Delta$ so that the pushoff $\partial F_u\subset \Delta$ is an embedding of $\partial F$ disjoint from $\partial F$. Extend $u$ to a normal vector field to $F$ that is transverse to $F$. Again we may assume that $F_u \cap \partial B$ is a section $s_0'$ of $E$ over an arc in the base sphere connecting the branch points and extend $s_0'$ to a section $\Sigma'$ of $E$ that intersects the zero section transversely in one point. The preimage of $F_u - \operatorname{int} B$ along with a half of the preimage of $\Sigma'$ then determines a pushoff $\tf'$ of $\tf$ in $\tw$. Note that any intersections between $F$ and $F_u$ contribute an even number to the intersection number $I=\tf \cdot \tf'$ thus the parity of $I$ depends on the intersections between the lifts of $\Sigma$ and $\Sigma'$. If $p$ is a positive double point, then $s_0'$ can be thought of as $s_0$ rotated by a small angle, so multiplied by a unit complex number $\xi$ close to 1. Hence we can choose $\Sigma'$ to be $\xi \Sigma$ and the contribution of the lifts of $\Sigma$ and $\Sigma'$ to $I$ is $\pm 1$. If on the other hand $p$ is a negative double point, then the pushoff $s_0'$ at the endpoints of the arc lies on the opposite sides of $s_0$ and this additional linking between the lifts of $s_0$ and $s_0'$ contributes another $\pm 1$ to $I$ making it even in this case.   
\endproof

The proof of Theorem \ref{thm:clasp} is now immediate upon recalling that by a standard algebraic topology argument (see for example \cite[Lemma 2.1]{qhs}) the determinant of the intersection form of a four-manifold that bounds a rational homology sphere divides the order of the first homology of the boundary (which for the double branched cover of a knot is equal to the determinant of the knot) with quotient a square.

\proof[Proof of Corollary \ref{cor:c=1}]
If $\sigma(K)=2$ and $K$ bounds an immersed disk in the four-ball with one double point, then it must be a negative double point.  This follows from Proposition \ref{prop:movie} and the fact that changing a positive crossing cannot decrease the signature (see \cite[Proposition 2.1]{cl-unknotting}, also \cite[Theorem 5.1]{st}).  
According to Theorem \ref{thm:clasp}, the double branched cover $Y$ of $K$ bounds a positive definite four-manifold $X$ with intersection form $Q=\left[\begin{matrix} a & 1\\ 1 & 2\end{matrix}\right]$ which presents a cyclic subgroup of $H_1(Y;\zz)$. As noted above, the determinant $r$ of $Q$ divides $\det K=|H_1(Y;\zz)|$ with quotient a square, $s^2$, and the image of the restriction homomorphism $\spinc(X) \to \spinc(Y)$ has order $rs$. For each \spinc structure $\spinct$ on $Y$ that extends over $X$ it follows from \cite[Theorems 1.2 and 9.6]{os4} that
\begin{equation}\label{eq:min}
\td(Y,\spinct):=-d(Y,\spinct)+\min_{\spincs}\frac{c_1(\spincs)^2-2}4
\end{equation}
is a nonnegative even integer, where the minimum is taken over all $\spincs\in\spinc(X)$ that restrict to $\spinct$.  Changing the \spinc structure on $X$ by a torsion element does not change the square of its Chern class.   The formula follows by noting that the minimum computes the $d$-invariant of a \spinc structure on $r/2$ surgery on the unknot \cite[Corollary 1.5]{plumbed}.  By the recursive formula for correction terms of lens spaces given in \cite{os4}, with the labelling shifted by $(r+1)/2$ so that the spin structure on the lens space is labelled by $i=0$, these are
\begin{equation}\label{eq:dlens}
d(S^3_{r/2}(O),i)=i^2/2r-\begin{cases}0& i \equiv (r-1)/2 \pmod 2\\1/2&i \equiv (r+1)/2 \pmod 2\end{cases}.
\end{equation}
(Alternatively, one can compute the minimum in \eqref{eq:min} directly using a suitable set of short characteristic covectors for the form $Q$ as in \cite{qhs} or \cite{osu1}.)

If $\sigma(K)=0$, then at least one of $K$ or its mirror $\overline{K}$ bounds a disk with a positive double point (and no other double points) to which the above may be applied.
\endproof


\section{Heegaard-Floer obstructions to (concordance) unknotting number one}
\label{sec:obstructions}

In this section we prove Theorems \ref{thm:u=1} and \ref{thm:uc=1} and Corollary \ref{cor:knots}.

\proof[Proof of Theorem \ref{thm:u=1}]
The first three conditions in the theorem are due to \ozsvath\ and \szabo\ \cite{osu1}. For completeness we give a proof of all five based on the formula of Ni and Wu \cite[Proposition 1.6]{niwu} for $d$-invariants of a positive Dehn surgery on a knot $K$ in the three-sphere: for $p,q>0$ 
\begin{equation}\label{eq:niwu}
d(S^3_{p/q}(K),i)=d(S^3_{p/q}(O),i)-2\max\{V_{\lfloor i/q \rfloor},H_{\lfloor (i-p)/q \rfloor}\}, \qquad i=0,\ldots,p-1,
\end{equation}
where $i$ enumerates relative \spinc structures on the complement of $K$ and $O$ denotes the unknot. We define normalized $d$-invariants as
\begin{equation}\label{eq:td}
\td_i:=-d(S^3_{p/q}(K),i)+d(S^3_{p/q}(O),i)=2\max\{V_{\lfloor i/q \rfloor},H_{\lfloor (i-p)/q \rfloor}\},
\end{equation}
for $i=0,\ldots,p-1$.

The sequences $V_j$ and $H_j$ ($j\in \zz$) are determined by the knot Floer chain complex $C=CFK^\infty(K)$ of $K$ and we recall their construction. In \cite{Zsurg} \ozsvath\ and \szabo\ define quotient complexes $A^+_j=C\{k\ge0\ \text{or}\ l\ge j\}$ and $B^+_j=B^+=C\{k\ge0\}$, where $(k,l)\in \zz \oplus \zz$ denotes the bidegree on $C$, and $U$-equivariant maps $v_j^+ \colon A_j^+ \to B^+$ and $h_j^+ \colon A_j^+ \to B^+$. The $U$-action lowers the bidegree by $(1,1)$ and $v_j^+$ is the obvious projection to the quotient complex; we will not need the description of $h_j^+$ but note that it is essentially the projection to $C\{l\ge j\}$. The homology of $B^+$ is isomorphic to $HF^+(S^3)\cong\calt^+=\zz[U,U^{-1}]/U\zz[U]$. Similarly there is a $\calt^+$ summand, which we denote $\calt_j^+$,  in the homology of $A_j^+$ (noting that the homology of $C$ is isomorphic to $HF^\infty(S^3)\cong \calt=\zz[U,U^{-1}]$, and the projection from $C$ to $A_j^+$ induces an isomorphism on homology in high degrees). 
Since both $v_j^+$ and $h_j^+$ induce isomorphisms on the chain level in sufficiently high degrees, their induced maps on homologies, restricted to $\calt^+_j$, are given by multiplication by $U^{V_j}$ and $U^{H_j}$. Moreover, for $j$ greater than the genus of $K$, $A_j^+=B^+$ and hence $V_j=0$ in this range. We restrict our analysis to $V_j$ for $j \ge 0$ as these determine $H_{-j}$ (see below) and all of these together are the only values of $V$'s and $H$'s relevant to \eqref{eq:niwu}. 

We first show that $H_{-j}=V_j$ for $j\ge 0$. (See also \cite[\S5.2]{niwuQgenus}.)
As noted in \cite[Corollary 2.3]{rs-links}, if $n$ is a large positive integer, then
$V_i=0$ for $i>n/2$ and $H_{i-n}=0$ for $i\le n/2$, and thus
the normalized $d$-invariants $\td_i$ of $S^3_n(K)$ are equal to either $2V_i$ for $i\le n/2$ or $2H_{i-n}$ for $i>n/2$. Since in the case of an integer surgery $i=0$ corresponds to a spin structure, conjugation invariance $\td_i=\td_{n-i}$ of (normalized) $d$-invariants yields the stated equality.

Next we claim that $V_j -1 \le V_{j+1} \le V_j$ for all $j$. Let $K_j$ denote the kernel of the projection $A_j^+ \to A_{j+1}^+$. Since the $U$-action on the left complex in the $U$-equivariant short exact sequence
$$ 0 \to K_j \to A_j^+ \to A_{j+1}^+ \to 0$$
is trivial, it follows that in homology the map from $\calt^+_{j+1}$ is trivial and the map into $\calt^+_j$ can contain in its image only the kernel of the $U$-action. Hence the map $\calt^+_j \to \calt^+_{j+1}$ is either the identity or multiplication by $U$. This proves the claim.

We now specialize to $p/q=r/2$. At this point we know that $V_j$ ($j\ge 0$) is a nonincreasing and $H_j$ ($j\le 0$) a nondecreasing sequence of nonnegative integers which yields the first two conditions. From the monotonicity properties of $V_j$ and $H_j$ and the fact that $V_0=H_0$ it follows that $\td_{2i}=\td_{2i+1}=2V_i$ and $\td_{r-2i}=\td_{r-2i+1}=2H_{-i}$ for small $i\ge 0$. Note that in this labelling of \spinc structures the spin structure is labelled by $i=(r+1)/2$. We claim that the symmetry condition above extends to the spin structure from both ends. Denote by $\ell$ the largest value of $i$ for which $V_{\lfloor i/2 \rfloor} \ge H_{\lfloor (i-r)/2 \rfloor}$. If $\ell \ge (r+1)/2$, then the symmetry condition holds up to the spin structure whereas for $\ell < (r+1)/2$ it holds above it. In both cases using conjugation invariance of $d$-invariants it follows that the symmetry holds on both sides of the spin structure. 

Next note that the monotonicity conditions on $V_j$ and $H_j$ imply that $\td_i$ is nonincreasing up to the spin structure (the symmetry shows us that $\td_{2i}=\td_{2i+1}$, and the inequality $\td_{2i+1}\le\td_{2i+2}$ in this range follows immediately from \eqref{eq:td}) and nondecreasing after that. Moreover, since $V_j$ and $H_j$ can jump by at most 1, the last condition also follows.

Finally we substitute for $d(S^3_{r/2}(O),i)$ using \eqref{eq:dlens}.
The isomorphism $\phi$ accounts for the fact that the given labelling of \spinc structures on $Y$ may not agree with the one assumed above.
\endproof

\proof[Proof of Theorem \ref{thm:uc=1}]
Let $Y'=S^3_{l/2}(C)$ with $l>0$ and odd, and let $W$ denote the rational homology cobordism with $\partial W=Y\sqcup-Y'$.  Let $X$ denote the surgery nullcobordism  bounded by $Y'$; this is given by a Kirby diagram with framing $(l+1)/2$ on $C$ and framing $2$ on a meridian of $C$.  The restriction map from $\spinc(X)$ to $\spinc(Y')$ is surjective since $H^2(X;\zz)$ is torsion-free.  We need to understand the image of the restriction map $\spinc(W)\to \spinc(Y) \times \spinc(Y')$, or equivalently from $\spinc(X\cup W)$. This image is conjugation invariant, and it has an odd number of elements, since it is affine isomorphic to a subgroup of $H^2(Y;\zz)\oplus H^2(Y';\zz)$.  Thus it contains a conjugation fixed element, which is the pair of spin structures on the two boundary components. It follows that it suffices to understand the restriction map on $H^2$.  Our goal is to understand the image of $H^2(X \cup W;\zz)$ in each of $H^2(Y';\zz)$ and $H^2(Y;\zz)$ separately; we will see in particular that the torsion subgroup of $H^2(X \cup W;\zz)$ restricts trivially to $Y'$ and may restrict nontrivially to $Y$.

Let $\Lambda_X=(H_2(X;\zz),Q_X)$ denote the intersection lattice of $X$, with dual lattice $\Lambda_X^*$ which is $H^2(X;\zz)$ with the induced pairing.
Let $\Lambda=(H_2(X \cup W;\zz)/\tors,Q_{X \cup W})$ denote the intersection lattice of the manifold $X\cup W$ bounded by $Y$, and let $\Lambda^*$ be its dual lattice (which is $H^2(X\cup W;\zz)/\tors$ with the induced pairing).  We have natural inclusions
$$\Lambda_X\subset\Lambda\subset\Lambda^*\subset\Lambda_X^*,$$
the last of which comes from the restriction map $H^2(X\cup W;\zz)\to H^2(X;\zz)$.
Let $r$ be the determinant of $\Lambda$ (that is, the order of $\Lambda^*/\Lambda$) and let $t$ be the index of $\Lambda_X$ in $\Lambda$.
Using the long exact sequence of the pair $(X,Y')$ we see that $\Lambda_X^*/\Lambda_X\cong H^2(Y';\zz)$ is cyclic of order $l=rt^2$.  The image of $H^2(X\cup W;\zz)$ in $H^2(Y';\zz)$ is the index $t$ subgroup $H'\cong\Lambda^*/\Lambda_X\subset \Lambda_X^*/\Lambda_X$.

Using the long exact sequence of the pair $(X\cup W,Y)$, as in for example \cite{qhs}, the order of {red}$H^2(Y;\zz)$ is $rs^2$ for some $s$.  Moreover, the image $H$ of $H^2(X\cup W;\zz)$ in $H^2(Y;\zz)$ has order $rs$, the image $T$ of the torsion subgroup of $H^2(X\cup W;\zz)$ in $H$ has order $s$, and $H/T \cong \Lambda^*/\Lambda\cong\zz/r\zz$.
The composite surjection $H^2(X\cup W;\zz)\to \Lambda^*\to\Lambda^*/\Lambda\cong\zz/r\zz$ thus factors through each of $H$ and $H'$.

Given a \spinc structure on $X\cup W$, its restriction to $Y$ may be changed by any given element of $T$ using the action of the torsion subgroup of $H^2(X\cup W;\zz)$, leaving the restriction to $Y'$ unchanged.  Since $d$-invariants are \spinc rational homology cobordism invariants \cite{os4}, this proves the first statement of the theorem.  Choose a labelling $\spinc(Y')\cong\zz/rt^2\zz$ as in Theorem \ref{thm:u=1}.
For each $i\in \zz/r\zz$ choose a \spinc structure $\spincs_i \in \spinc(X\cup W)$ whose restriction to $Y'$ is labelled by $ti$, and let $\phi(\spincs_i|_Y)=i$.
From the discussion above, this is the quotient homomorphism $H\to H/T$ composed with an automorphism of $\zz/r\zz$. 
Then 
$$d(Y',ti)=d(Y,\spincs_i|_Y)$$
and the second conclusion of the theorem now follows from the conclusion of Theorem \ref{thm:u=1} applied to $Y'$.
\endproof

Note that since in the above proof we are using the conclusions of Theorem \ref{thm:u=1} on a subgroup of \spinc structures on the surgery manifold, we lose the symmetry and boundedness conditions from Theorem \ref{thm:u=1}.

In the proof of Corollary \ref{cor:knots} we need to know whether the double branched cover of an unknotting number one knot is a positive or a negative surgery on some knot. The exact information is given by the following signed refinement of Montesinos' trick, which is proved in \cite[Theorem 8.1]{osu1}. However, we only need to pin down the sign in case the signature is nonzero and we give an alternative argument for that. 
\begin{lemma}\label{lem:u=1}
Let $K$ be a knot in $S^3$ with determinant $l$ and unknotting number one.  Suppose that either
\begin{enumerate}[(i)]
\item $K$ has signature two, or
\item $K$ has signature zero and may be unknotted by changing a positive crossing to a negative crossing.
\end{enumerate}
Then $\Sigma(K)=S^3_{l/2}(C)$.
\end{lemma}

\proof[Proof of (i)]
If the signature of $K$ is two, then $Y=\Sigma(K)$ is the boundary of a spin four-manifold $W$ with signature two. Indeed, $W$ can either be taken to be the double branched cover of a Seifert surface for $K$ pushed into the four-ball or the manifold from Theorem \ref{thm:clasp}. In the latter case note that since by \cite[Theorem 3.7]{cl-unknotting} $H_1(W;\zz/2\zz)=0$ and the intersection form is even (only changing a negative crossing can unknot $K$), the manifold is spin. Moreover, since $\sigma(K)=2$, $\det(K) \equiv -1 \pmod 4$ by \cite[Theorem 5.6]{murasugi}, the surgery manifold $X$, corresponding to the positive surgery, that $\pm Y$ bounds, is also spin with signature two. If $-Y$ were the positive surgery, then $W \cup X$ would be a 
closed spin manifold of signature four, a contradiction.
\endproof

\proof[Proof of Corollary \ref{cor:knots}]
Suppose $K$ is a knot with unknotting number one and let $Y=\Sigma(K)$. If $K$ has signature zero, then either $Y$ or $\Sigma(\overline{K})=-Y$ is a positive surgery on a knot; if necessary we replace $K$ by its mirror so that $Y$ is a positive surgery. If $K$ has signature two, then it follows from Lemma \ref{lem:u=1}(i) that $Y$ is a positive surgery. 
In each case we then apply Theorem \ref{thm:u=1}.

Now suppose $K$ is concordant to a knot $K'$ with unknotting number one.
The argument in the preceding paragraph may be applied to $K'$, which has the same signature as $K$, and we then apply Theorem \ref{thm:uc=1}, using the fact that the double branched cover $W$ of a concordance between $K$ and $K'$ gives a rational homology cobordism between $\Sigma(K)$ and $\Sigma(K')$.
\endproof


\section{Low crossing number examples}
\label{sec:examples}

Table \ref{table:smallknots} lists the determinant, unknotting number, concordance unknotting number, slicing number, four-ball crossing number, four-ball genus and half of the absolute value of the signature for knots with up to 9 crossings, with incomplete information for three knots.  In this section we describe how the previously unknown values in the table have been obtained.

The slicing number of $7_4$ was determined by Livingston \cite{liv} and the four-ball crossing number of $8_{16}$ was determined by Murakami and Yasuhara \cite{my}.  Apart from the values of $u_c$, $u_s$ and $c^*$, the data in the table is taken from \cite{knotinfo}.

For any knot for which the unknotting number is equal to the four-ball genus, it follows from \eqref{eq:ineq} that $u=u_c=u_s=c^*=g^*$.  Also any slice knot has $u_c=u_s=c^*=g^*=0$.  We have highlighted the knots in Table \ref{table:smallknots} for which neither of these situations apply, and we give some details about each in what follows.  

At the end of the section we give a brief account of the computation of the invariants in Table \ref{table:10cross} of knots with 10 crossings.

\subsection{Calculation of $d$-invariants}  For the knots we consider in this section, we calculate the correction terms of the double branched cover using various methods from \cite{os4} and \cite{osu1}.  The recursive formula in \cite[Section 4.1]{os4}
may be used to compute $d$-invariants of lens spaces, which are the double branched covers of 2-bridge knots and links.
More generally, we compute $d$-invariants using the intersection forms of sharp four-manifolds, as we recall now.

Let $Y$ be a rational homology 3-sphere, and suppose $Y$ (smoothly) bounds a positive-definite four-manifold $X$ with intersection form $Q$.  Then it follows from \cite[Theorems 1.2, 9.6]{os4} that for each $\spincs\in\spinc(X)$, we have
\begin{align}
\frac{c_1(\spincs)^2-b_2(X)}4&\ge d(Y,\spincs|_Y)\label{eq:thm9.6}\\
\mbox{and}\quad\frac{c_1(\spincs)^2-b_2(X)}4&\equiv d(Y,\spincs|_Y)\pmod2.\label{eq:thm1.2}
\end{align}

We say that $X$ is \emph{sharp} if every \spinc structure on $Y=\partial X$ admits an extension over $X$ with equality in \eqref{eq:thm9.6}; in this case we may use the intersection form on $X$ to calculate the correction terms of $Y$.  When $Y$ is the double branched cover of an alternating knot $K$, it is shown in \cite{osu1} that the Goeritz matrix of an alternating diagram of $K$ is the intersection form of a sharp manifold bounded by $Y$.  A method is given in the same paper to exhibit sharp four-manifolds bounded by the double branched covers of certain nonalternating knots via exact triangles; this applies in particular to $K=10_{158}$ which we will consider in Subsection \ref{subsec:10}.

Suppose then that a positive-definite symmetric integer matrix $G=(g_{ij})$ of rank $n$ is the intersection form of a sharp four-manifold $X$ bounded by $Y$.
Call $\xi=(\xi_1,\dots,\xi_n)$ a characteristic covector for $G$ if $$\xi_i\equiv g_{ii}\pmod2, i=1,\dots,n,$$
and denote the set of all characteristic covectors by $\Char(G)$; these represent first Chern classes of \spinc structures on $X$.  The set $\Char(G)/2G\zz^n$ is then affine isomorphic to $\spinc(Y)$.
Using a computer, we  partition the set of characteristic covectors $\xi=(\xi_1,\dots,\xi_n)$ with $g_{ii}\le \xi_i<-g_{ii}$ into cosets of $\Char(G)/2G\zz^n$, and use the Smith normal form of $G$ to record the affine group structure of $\Char(G)/2G\zz^n$.  We also minimise the quantity 
$(\xi^TQ^{-1}\xi-n)/4$ on each of these finite sets of coset representatives; the resulting rational numbers associated to $\spinc(Y)$ are the $d$-invariants.

\subsection{Application of Theorem \ref{thm:clasp}.}
In applications, one wishes to determine whether a given rational homology sphere $Y$ may bound a positive-definite 4-manifold whose intersection form is of half-integer surgery type.  One can obstruct a given form using the $d$-invariants of $Y$, together with \eqref{eq:thm9.6} and \eqref{eq:thm1.2}.  There are finitely many positive-definite forms of a given rank and determinant, and a complete list of representatives may be found using the theory of reduced forms, see for example \cite{cassels,jones}.  For more details on how to write down a complete list of positive-definite forms of half-integer type, with a given rank and determinant, see \cite{unknotting}.

\subsection{Knots with $u=c^*$}
If a knot has $u=c^*$, then these also are equal to $u_c$ and $u_s$ by \eqref{eq:ineq}.
Using Theorem \ref{thm:clasp} in place of \cite[Theorem 2]{slicing} we may adapt Corollaries 3 and 4 of \cite{slicing} to conclude that
$7_4$, $8_{16}$, $9_5$, $9_{15}$, $9_{17}$ and $9_{31}$ all have $u=c^*=2$ and  $9_{10}$, $9_{13}$, and $9_{38}$ have $u=c^*=3$.  For the reader's convenience we will recall how the argument goes for one example.  Consider $K=9_{10}$.  This is a two-bridge knot with determinant $33$ and signature four.  The correction terms of the double branched cover $Y=L(33,23)$ comprise a function
$$d_Y:\zz/33\zz\to\qq,$$
well-defined up to a group automorphism of the domain, which may be computed using the recursion formula from \cite{os4}.  The number of negative double points
in a normally immersed disk in the four-ball bounded by $K$ is bounded below by half the signature (this follows for example from Theorem \ref{thm:imm}).  We wish to show that there does not exist such a disk with 2 double points; since the signature is four, these would both have to be negative and we may apply Theorem \ref{thm:clasp}.  Existence of such a disk then implies that $Y$ bounds a positive-definite four-manifold $X$ with $b_2=4$ and with an even intersection form of half-integer surgery type.  The determinant of this form has to be $33$ since by the long exact cohomology sequence of the pair $(X,Y)$ it divides the order of $H^2(Y;\zz)$ with quotient a square.  There turn out to be exactly two such forms, which are represented by the matrices
$$\left(\begin{matrix}
6 & 0 & 1 & 0\\
0 & 2 & 0 & 1\\
1 & 0 & 2 & 0\\
0 & 1 & 0 & 2
\end{matrix}\right),\quad
\left(\begin{matrix}
4 & 2 & 1 & 0\\
2 & 4 & 0 & 1\\
1 & 0 & 2 & 0\\
0 & 1 & 0 & 2
\end{matrix}\right).$$

Letting $\mathrm{Char}(Q)$ denote the set of elements of $\zz^4$ with even coefficients, the quotient of $\mathrm{Char}(Q)$ by the image of $2Q$ for each matrix $Q$ above is isomorphic to $\zz/33\zz$.   We may then compute the function
\begin{align*}
m_Q:\zz/33\zz&\to\qq\\
j &\mapsto \min\left\{\left.\frac{\xi^T Q^{-1}\xi-4}4\,\,\right|
\,\,\xi\in \mathrm{Char}(Q),\,[\xi]=j\right\}.
\end{align*}
(It suffices to consider vectors $\xi$ with entries bounded in absolute value by the corresponding diagonal element of $Q$, and minimise over the finite equivalence classes of such vectors.)
Using \cite{os4} we have that if $Y$ bounds $X$ with intersection form $Q$, then there exists an automorphism $\phi$ of $\zz/33\zz$ such that for each $j\in\zz/33\zz$, the difference
$$m_Q(j)-d_Y(\phi(j))$$ is a nonnegative even integer.  From inspection of $d_Y$ and each of the $m_Q$ it is straightforward to see that no such $\phi$ exists.  Indeed the correction terms of $Y$ are
$$\left\{
\begin{array}{rrrrrrrrrrr}
-1,& -\frac{23}{33},& \frac{7}{33},& -\frac{3}{11},& -\frac{5}{33},& \frac{19}{33},
& -\frac{1}{11},& -\frac{5}{33},& \frac{13}{33},& -\frac{5}{11},& -\frac{23}{33},\\
&&&&&&&&&&\\
-\frac{1}{3},& \frac{7}{11},& \frac{7}{33},& \frac{13}{33},& \frac{13}{11},
& \frac{19}{33},& \frac{19}{33},& \frac{13}{11},& \frac{13}{33},& \frac{7}{33},
& \frac{7}{11},\\
&&&&&&&&&&\\
-\frac{1}{3},& -\frac{23}{33},& -\frac{5}{11},& \frac{13}{33},& -\frac{5}{33},
& -\frac{1}{11},& \frac{19}{33},& -\frac{5}{33},& -\frac{3}{11},& \frac{7}{33},
& -\frac{23}{33}\\
\end{array}
\right\};$$
the order of the list corresponds to the group structure, but we will not need this for this example.

For the first matrix above we find there exists $j$ with $m_Q(j)=-9/11$, and for the second there exists $j$ with $m_Q(j)=-7/11$.  None of the correction terms of $Y$ have the property that $-9/11-d_Y(j)$ or $-7/11-d_Y(j)$ is a nonnegative even integer.  We conclude that $K$ has $c^*\ge3$, and in fact since $K$ can be unknotted by three crossing changes in its alternating diagram, that $u(K)=c^*(K)=3$.

For the knot $K=8_{18}$, the homology group of $Y=\Sigma(K)$ is isomorphic to $\zz/3\zz\oplus\zz/15\zz$.  The signature of this knot is zero.  
The only factorisation $45=rs^2$ for which there is an epimorphism from an order $rs$ subgroup of $\zz/3\zz\oplus\zz/15\zz$ onto $\zz/r\zz$ is $r=5$ and $s=3$.
The $d$-invariants of the unique order 5 subgroup of $\spinc(Y)$ may be calculated as in \cite[Proposition 3.2]{osu1} to be $$[0,4/5,-4/5,-4/5,4/5].$$
This subgroup gives a full set of representatives for the order 5 quotient of any of the possible order 15 subgroups of $\spinc(Y)$.
Thus by Corollary \ref{cor:c=1}, if $c^*(K)=1$ there would exist a positive even matching between these numbers (up to overall sign change) and the $d$-invariants of $5/2$ surgery on the unknot, which are
$$[0,-2/5,2/5,2/5,-2/5].$$
No such matching exists; in fact no integral matching exists, so that $Y$ is obstructed by the linking form from bounding a \emph{topological} manifold with intersection form $Q=\left(\begin{matrix}3&1\\1&2\\\end{matrix}\right)$.  We conclude that $u=c^*=2$ for this knot.

Similar reasoning applies to the knot $9_{40}$ for which the homology group is $\zz/5\zz\oplus\zz/15\zz$ and the $d$-invariants on the unique $\zz/3\zz$ subgroup are $[-1/2,5/6,5/6]$, while those of $3/2$ surgery on the unknot are $[-1/2,1/6,1/6]$.

The knot $K=9_{49}$ has signature four and determinant 25.  If $c^*(K)=2$, then as in the proof of Corollary \ref{cor:c=1}, $K$ must bound an immersed disk with $2$ negative (and no positive) double points.
Now by Theorem \ref{thm:clasp}, $\Sigma(K)$ bounds a positive-definite even form of half-integer surgery type of rank 4 and determinant 1 or 25.  However no such forms exist (this may be checked as in \cite[\S6]{unknotting}) and we conclude that $u=c^*=3$.

\subsection{Knots with $u_c=u_s=c^*=1$}
The knots $8_{10}$ and $9_{37}$ each have $u_c=u_s=1$ as may be seen from Figure \ref{fig:knots}: the former is concordant to the trefoil and is one crossing change from $3_1\#\overline{3_1}$, the latter is concordant to $4_1$ and gives $6_1$ after changing the indicated crossing.  Comparing with the proof of Proposition \ref{prop:movie} one may observe that in each case the crossing change and the concordance combine to give an immersed disk with a single double point.

\subsection{Knots with $u_c=2$, $u_s=c^*=1$}
The knots listed in Table \ref{table:data} each have $u_s=1$ as may be seen by changing the crossing circled in Figure \ref{fig:knots}; in each case the slice knot thus obtained is given in the table.

\begin{table}[!ht]
\begin{center}
\begin{tabular}{||c|c|c||}
\hline
\rule{0pt}{3ex} 
Knot  & $\td$ & Slice knot\\
\hline
\rule{0pt}{2.5ex} 
$8_3$ & $0^6,2,0,2$ & $6_1$ \\
\hline
\rule{0pt}{2.5ex} 
$8_4$ & $0^5,2,0,2^3$ & $6_1$ \\
\hline
\rule{0pt}{2.5ex} 
$8_6$ & $0^7,2,0,2^3$ & $6_1$ \\
\hline
\rule{0pt}{2.5ex} 
$8_{12}$ & $0^8,2,0,2^4,4$ & $6_1$ \\
\hline
\rule{0pt}{2.5ex} 
$9_8$ & $0^7,2,0,2^4,4,2,4$ & $6_1$ \\
\hline
\rule{0pt}{2.5ex} 
$9_{25}$ & $0^{11},2,0,2^4,4,2,4^4,6$ & $6_1$ \\
\hline
\rule{0pt}{2.5ex} 
$9_{29}$ & $0^9,2,0,2^6,4,2,4^4,6^3$ & $3_1\#\overline{3_1}$ \\
\hline
\rule{0pt}{2.5ex} 
$9_{32}$ & $0^{11},2,0,2^6,4,2,4^4,6^4,8$ & $6_1$ \\
\hline
\end{tabular}
\vskip5mm
\begin{narrow}{0.3in}{0.3in}
\caption{
{\bf{Data for knots with $u_c=2$ and $u_s=1$, and with cyclic $H_1(\Sigma(K))$.}}
In the second column of the table we have used an abbreviated notation in which for example $0^6$ stands for $0,0,0,0,0,0$.
}
\label{table:data}
\end{narrow}
\end{center}
\end{table}

We consider first the knot $K=8_3$, which is the two-bridge knot $S(17,4)$ with determinant $17$ and signature zero.  The $d$-invariants may be computed using either the recursive formula from \cite{os4} or the Goeritz matrix from an alternating diagram as in \cite{osu1}.  In cyclic group order starting at the spin structure these are
$$\textstyle[0,\frac{4}{17},\frac{16}{17},\frac{2}{17},-\frac{4}{17},-\frac{2}{17},\frac{8}{17},-\frac{8}{17},-\frac{16}{17},-\frac{16}{17},-\frac{8}{17},\frac{8}{17},-\frac{2}{17},-\frac{4}{17},\frac{2}{17},\frac{16}{17},\frac{4}{17}].$$
Since the determinant is square-free we are looking for a positive even matching on the whole group; we must compare these $d$-invariants with those of $17/2$ surgery on the unknot which are
$$\textstyle\left[0,-\frac{8}{17},\frac{2}{17},-\frac{4}{17},\frac{8}{17},\frac{4}{17},\frac{18}{17},\frac{16}{17},\frac{32}{17},\frac{32}{17},\frac{16}{17},\frac{18}{17},\frac{4}{17},\frac{8}{17},-\frac{4}{17},\frac{2}{17},-\frac{8}{17}\right].$$
Since the signature is zero we are free to switch the sign of the first list above.  We find there are two positive even matchings $\phi:\zz/17\zz\to\zz/17\zz$.  The first is multiplication by 5 and applies to the list above, and the second is multiplication by 3 and applies to the list with the opposite sign.  Both result in the same ordered list of $\td$ invariants as in Theorem \ref{thm:uc=1} which are
$$0,0,0,0,0,0,2,0,2.$$
This fails monotonicity and so $u_c(K)=u(K)=2$.

Similar analysis applies to each of the knots listed in Table \ref{table:data}.

\begin{table}[!ht]
\tiny
\begin{center}
\begin{tabular}{||c|c|c|c|c|c|c|c||c|c|c|c|c|c|c|c||}
\hline
\rule{0pt}{4ex} 
\rule[-2.3ex]{0pt}{0pt}
Knot  & $\det$ & $u$ & $u_c$ & $u_s$ & $c^*$ & $g^*$ & $\dfrac{|\sigma|}2$ & Knot & $\det$ & $u$ & $u_c$ & $u_s$ & $c^*$ & $g^*$ & $\dfrac{|\sigma|}2$\\
\hline
$3_1$ & $3$ & $1$ & $1$ & $1$ & $1$ & $1$ & $1$ & \gc$9_8$ & $31$ & $2$ & $2$ & $1$ & $1$ & $1$ & $1$ \\
\hline
$4_1$ & $5$ & $1$ & $1$ & $1$ & $1$ & $1$ & $0$ & $9_9$ & $31$ & $3$ & $3$ & $3$ & $3$ & $3$ & $3$ \\
\hline
$5_1$ & $5$ & $2$ & $2$ & $2$ & $2$ & $2$ & $2$ & \gc$9_{10}$ & $33$ & $3$ & $3$ & $3$ & $3$ & $2$ & $2$ \\
\hline
$5_2$ & $7$ & $1$ & $1$ & $1$ & $1$ & $1$ & $1$ & $9_{11}$ & $33$ & $2$ & $2$ & $2$ & $2$ & $2$ & $2$ \\
\hline
$6_1$ & $9$ & $1$ & $0$ & $0$ & $0$ & $0$ & $0$ & $9_{12}$ & $35$ & $1$ & $1$ & $1$ & $1$ & $1$ & $1$ \\
\hline
$6_2$ & $11$ & $1$ & $1$ & $1$ & $1$ & $1$ & $1$ & \gc$9_{13}$ & $37$ & $3$ & $3$ & $3$ & $3$ & $2$ & $2$ \\
\hline
$6_3$ & $13$ & $1$ & $1$ & $1$ & $1$ & $1$ & $0$ & $9_{14}$ & $37$ & $1$ & $1$ & $1$ & $1$ & $1$ & $0$ \\
\hline
$7_1$ & $7$ & $3$ & $3$ & $3$ & $3$ & $3$ & $3$ & \gc$9_{15}$ & $39$ & $2$ & $2$ & $2$ & $2$ & $1$ & $1$ \\
\hline
$7_2$ & $11$ & $1$ & $1$ & $1$ & $1$ & $1$ & $1$ & $9_{16}$ & $39$ & $3$ & $3$ & $3$ & $3$ & $3$ & $3$ \\
\hline
$7_3$ & $13$ & $2$ & $2$ & $2$ & $2$ & $2$ & $2$ & \gc$9_{17}$ & $39$ & $2$ & $2$ & $2$ & $2$ & $1$ & $1$ \\
\hline
\gc$7_4$ & $15$ & $2$ & $2$ & $2$ & $2$ & $1$ & $1$ & $9_{18}$ & $41$ & $2$ & $2$ & $2$ & $2$ & $2$ & $2$ \\
\hline
$7_5$ & $17$ & $2$ & $2$ & $2$ & $2$ & $2$ & $2$ & $9_{19}$ & $41$ & $1$ & $1$ & $1$ & $1$ & $1$ & $0$ \\
\hline
$7_6$ & $19$ & $1$ & $1$ & $1$ & $1$ & $1$ & $1$ & $9_{20}$ & $41$ & $2$ & $2$ & $2$ & $2$ & $2$ & $2$ \\
\hline
$7_7$ & $21$ & $1$ & $1$ & $1$ & $1$ & $1$ & $0$ & $9_{21}$ & $43$ & $1$ & $1$ & $1$ & $1$ & $1$ & $1$ \\
\hline
$8_1$ & $13$ & $1$ & $1$ & $1$ & $1$ & $1$ & $0$ & $9_{22}$ & $43$ & $1$ & $1$ & $1$ & $1$ & $1$ & $1$ \\
\hline
$8_2$ & $17$ & $2$ & $2$ & $2$ & $2$ & $2$ & $2$ & $9_{23}$ & $45$ & $2$ & $2$ & $2$ & $2$ & $2$ & $2$ \\
\hline
\gc$8_3$ & $17$ & $2$ & $2$ & $1$ & $1$ & $1$ & $0$ & $9_{24}$ & $45$ & $1$ & $1$ & $1$ & $1$ & $1$ & $0$ \\
\hline
\gc$8_4$ & $19$ & $2$ & $2$ & $1$ & $1$ & $1$ & $1$ & \gc$9_{25}$ & $47$ & $2$ & $2$ & $1$ & $1$ & $1$ & $1$ \\
\hline
$8_5$ & $21$ & $2$ & $2$ & $2$ & $2$ & $2$ & $2$ & $9_{26}$ & $47$ & $1$ & $1$ & $1$ & $1$ & $1$ & $1$ \\
\hline
\gc$8_6$ & $23$ & $2$ & $2$ & $1$ & $1$ & $1$ & $1$ & $9_{27}$ & $49$ & $1$ & $0$ & $0$ & $0$ & $0$ & $0$ \\
\hline
$8_7$ & $23$ & $1$ & $1$ & $1$ & $1$ & $1$ & $1$ & $9_{28}$ & $51$ & $1$ & $1$ & $1$ & $1$ & $1$ & $1$ \\
\hline
$8_8$ & $25$ & $2$ & $0$ & $0$ & $0$ & $0$ & $0$ & \gc$9_{29}$ & $51$ & $2$ & $2$ & $1$ & $1$ & $1$ & $1$ \\
\hline
$8_9$ & $25$ & $1$ & $0$ & $0$ & $0$ & $0$ & $0$ & $9_{30}$ & $53$ & $1$ & $1$ & $1$ & $1$ & $1$ & $0$ \\
\hline
\gc$8_{10}$ & $27$ & $2$ & $1$ & $1$ & $1$ & $1$ & $1$ & \gc$9_{31}$ & $55$ & $2$ & $2$ & $2$ & $2$ & $1$ & $1$ \\
\hline
$8_{11}$ & $27$ & $1$ & $1$ & $1$ & $1$ & $1$ & $1$ & \gc$9_{32}$ & $59$ & $2$ & $2$ & $1$ & $1$ & $1$ & $1$ \\
\hline
\gc$8_{12}$ & $29$ & $2$ & $2$ & $1$ & $1$ & $1$ & $0$ & $9_{33}$ & $61$ & $1$ & $1$ & $1$ & $1$ & $1$ & $0$ \\
\hline
$8_{13}$ & $29$ & $1$ & $1$ & $1$ & $1$ & $1$ & $0$ & $9_{34}$ & $69$ & $1$ & $1$ & $1$ & $1$ & $1$ & $0$ \\
\hline
$8_{14}$ & $31$ & $1$ & $1$ & $1$ & $1$ & $1$ & $1$ & \gc$9_{35}$ & $27$ & $3$ & $[2,3]$ & $2$ & $2$ & $1$ & $1$ \\
\hline
$8_{15}$ & $33$ & $2$ & $2$ & $2$ & $2$ & $2$ & $2$ & $9_{36}$ & $37$ & $2$ & $2$ & $2$ & $2$ & $2$ & $2$ \\
\hline
\gc$8_{16}$ & $35$ & $2$ & $2$ & $2$ & $2$ & $1$ & $1$ & \gc$9_{37}$ & $45$ & $2$ & $1$ & $1$ & $1$ & $1$ & $0$ \\
\hline
$8_{17}$ & $37$ & $1$ & $1$ & $1$ & $1$ & $1$ & $0$ & \gc$9_{38}$ & $57$ & $3$ & $3$ & $3$ & $3$ & $2$ & $2$ \\
\hline
\gc$8_{18}$ & $45$ & $2$ & $2$ & $2$ & $2$ & $1$ & $0$ & $9_{39}$ & $55$ & $1$ & $1$ & $1$ & $1$ & $1$ & $1$ \\
\hline
$8_{19}$ & $3$ & $3$ & $3$ & $3$ & $3$ & $3$ & $3$ & \gc$9_{40}$ & $75$ & $2$ & $2$ & $2$ & $2$ & $1$ & $1$ \\
\hline
$8_{20}$ & $9$ & $1$ & $0$ & $0$ & $0$ & $0$ & $0$ & $9_{41}$ & $49$ & $2$ & $0$ & $0$ & $0$ & $0$ & $0$ \\
\hline
$8_{21}$ & $15$ & $1$ & $1$ & $1$ & $1$ & $1$ & $1$ & $9_{42}$ & $7$ & $1$ & $1$ & $1$ & $1$ & $1$ & $1$ \\
\hline
$9_1$ & $9$ & $4$ & $4$ & $4$ & $4$ & $4$ & $4$ & $9_{43}$ & $13$ & $2$ & $2$ & $2$ & $2$ & $2$ & $2$ \\
\hline
$9_2$ & $15$ & $1$ & $1$ & $1$ & $1$ & $1$ & $1$ & $9_{44}$ & $17$ & $1$ & $1$ & $1$ & $1$ & $1$ & $0$ \\
\hline
$9_3$ & $19$ & $3$ & $3$ & $3$ & $3$ & $3$ & $3$ & $9_{45}$ & $23$ & $1$ & $1$ & $1$ & $1$ & $1$ & $1$ \\
\hline
$9_4$ & $21$ & $2$ & $2$ & $2$ & $2$ & $2$ & $2$ & $9_{46}$ & $9$ & $2$ & $0$ & $0$ & $0$ & $0$ & $0$ \\
\hline
\gc$9_5$ & $23$ & $2$ & $2$ & $2$ & $2$ & $1$ & $1$ & \gc$9_{47}$ & $27$ & $2$ & $[1,2]$ & $1$ & $1$ & $1$ & $1$ \\
\hline
$9_6$ & $27$ & $3$ & $3$ & $3$ & $3$ & $3$ & $3$ & \gc$9_{48}$ & $27$ & $2$ & $[1,2]$ & $1$ & $1$ & $1$ & $1$ \\
\hline
$9_7$ & $29$ & $2$ & $2$ & $2$ & $2$ & $2$ & $2$ & \gc$9_{49}$ & $25$ & $3$ & $3$ & $3$ & $3$ & $2$ & $2$ \\
\hline
\end{tabular}
\vskip5mm
\begin{narrow}{0.3in}{0.3in}
\caption{
{\bf Invariants of knots with at most 9 crossings.}  Knots for which calculation of $u_s$, $u_c$ or $c^*$ is nontrivial are highlighted.}
\label{table:smallknots}
\end{narrow}
\end{center}
\end{table}

\subsection{Knots with unknown values}
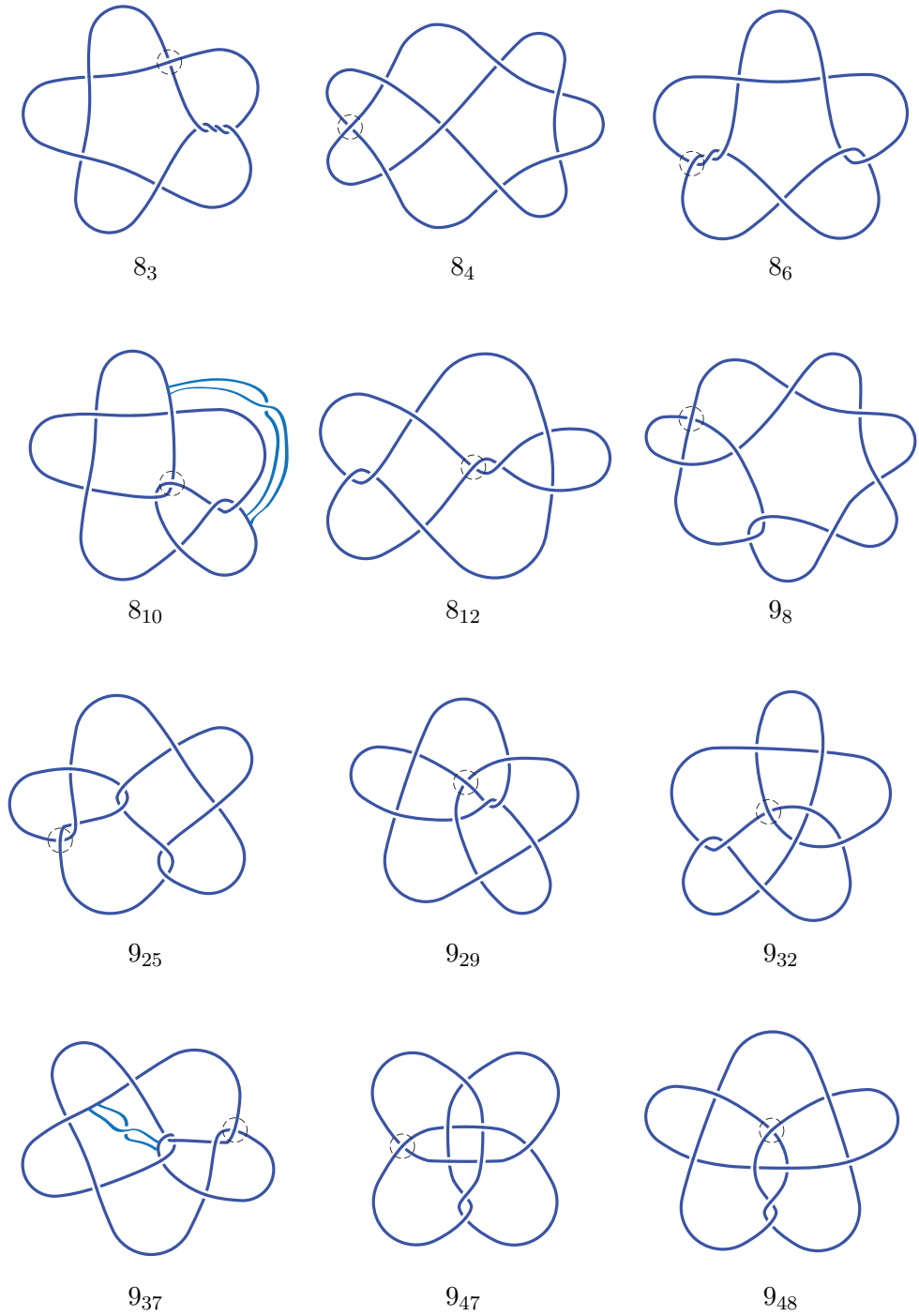
\begin{figure}[htbp]
\begin{center}
\ifpic{
\labellist
\small\hair 3pt
\pinlabel $8_3$ at 100 610
\pinlabel $8_4$ at 285 610
\pinlabel $8_6$ at 470 610
\pinlabel $8_{10}$ at 100 410
\pinlabel $8_{12}$ at 285 410
\pinlabel $9_8$ at 470 410
\pinlabel $9_{25}$ at 100 210
\pinlabel $9_{29}$ at 285 210
\pinlabel $9_{32}$ at 470 210
\pinlabel $9_{37}$ at 100 10
\pinlabel $9_{47}$ at 285 10
\pinlabel $9_{48}$ at 470 10
\endlabellist
\includegraphics[scale=0.7]{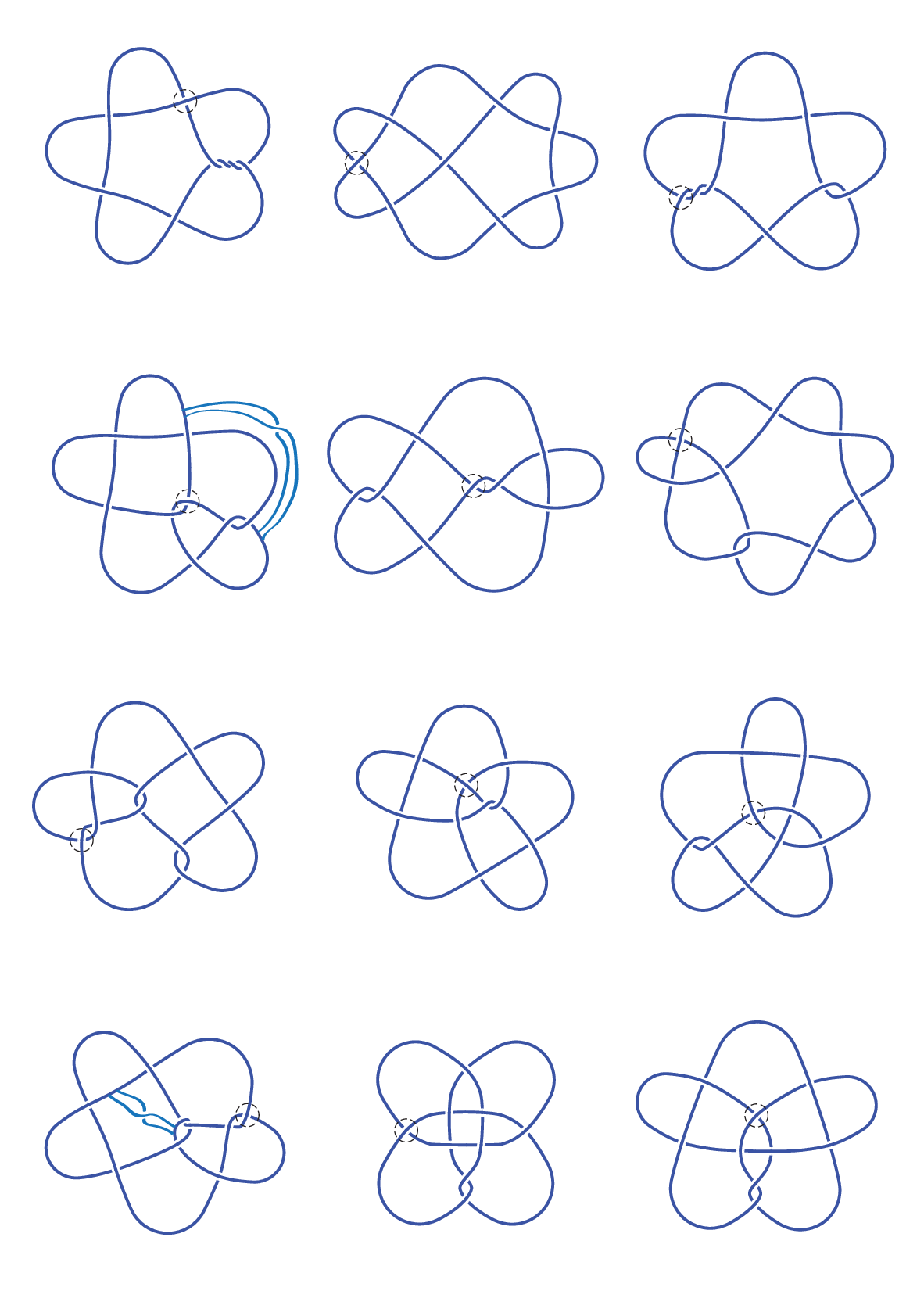}
}
\else \vskip 5cm \fi
\begin{narrow}{0.3in}{0.3in}
\caption
{\bf Crossing changes giving slice knots, and band moves giving concordances to unknotting number one knots.}
\label{fig:knots}
\end{narrow}
\end{center}
\end{figure}

The knots $9_{47}$ and $9_{48}$ have $u_s=1$ as may be seen by changing a crossing as in Figure \ref{fig:knots}, resulting in $6_1$ and $8_{20}$ respectively.

We next consider $K=9_{35}$, which is the pretzel knot $P(3,3,3)$.  This has signature two and $H_1(\Sigma(K))\cong\zz/3\zz\oplus\zz/9\zz$.   We give two proofs that $u_c(K)>1$.  The $d$-invariants (multiplied by 18) are
$$\left[
\begin{matrix}
-9 & 19 & -5 & 27 & 7 & 7 & 27 & -5 & 19\\
3 & -5 & 7 & 3 & 19 & 19 & 3 & 7 & -5\\
3 & -5 & 7 & 3 & 19 & 19 & 3 & 7 & -5
\end{matrix}
\right].$$
Here the rectangular array shows the group structure with the invariant coming from the spin structure in the top left position.  There are four subgroups of order 9; none of these admits an epimorphism onto $\zz/3\zz$ with constant $d$-invariant on fibres.  This is easy to see since the $d$-invariant of the spin structure is not repeated in any other \spinc structure.  Thus there does not exist a positive even subgroup matching, so by Corollary \ref{cor:knots} we see that $u_c(K)>1$.

We may also use a result from \cite{slicing} based on Donaldson's diagonalisation theorem \cite{donaldson} to show that in fact $c^*(K)>1$.
Combining the proof of \cite[Corollary 5]{slicing} with Theorem \ref{thm:clasp} we find that any normally immersed disk in $B^4$ bounded by the knot $7_4$ has at least two negative double points.\footnote{There is an oversight in the proof of \cite[Corollary 5]{slicing}: using the notation therein, there is more than one embedding of the lattice $L'_n$ in $\zz^m$, however the conclusion that the orthogonal complement does not contain a finite index sublattice of half-integer type is correct for any such embedding.}
Since $7_4$ is the $P(1,3,3)$ pretzel, we obtain $K=9_{35}$ from it by changing a positive crossing.  Thus an immersed disk bounded by $K$ gives rise to one for $7_4$ with one additional positive double point.  We conclude that $c^*(K)>1$; moreover, if $c^*(K)=2$, then any immersed disk realising this bound has two negative double points.  An argument of Traczyk \cite{traczyk} using the Jones polynomial shows that $9_{35}$ cannot be \emph{unknotted} by changing two negative crossings.  However as we see from Figure \ref{fig:moreknots} it is possible to go from $9_{35}$ to the slice knot $8_{20}$, which incidentally is the pretzel $P(3,2,-3)$, by two crossing changes.  Thus $u_s=c^*=2$ for this knot.

\subsection{Knots with 10 crossings.}\label{subsec:10}
Table \ref{table:10cross} lists invariants for 10 crossing knots; we have omitted slice knots or knots for which $u_c$, $u_s$ and $c^*$ are computable from \eqref{eq:ineq}.  Here we briefly indicate how the data in the table was compiled.

The knots $10_{19}$, $ 
10_{20}$, $ 10_{24}$, $ 10_{36}$, $  
10_{68}$, $ 10_{69}$, $ 10_{86}$, $10_{97}$, $ 10_{105}$, $ 
10_{109}$, $ 10_{116}$, $ 10_{121}$, $ 10_{122}$, $ 
10_{144}$, $ 10_{163}$, $ 10_{165}$
are obstructed  from having $c^*=1$; $10_{53}$, $10_{101}$ and $10_{120}$ are similarly obstructed  from having $c^*=2$.  This follows as in Corollaries 2 and 3 of \cite{slicing}, using Theorem \ref{thm:clasp} in place of \cite[Theorem 2]{slicing}.

The knots $10_{40}$, $10_{65}$, $10_{67}$, $10_{74}$, $10_{77}$, $10_{103}$ and $10_{106}$ are concordant to unknotting number one knots \cite{liv2} and hence have $u_c=1$. 

All remaining knots in Table \ref{table:10cross} with slice genus one are obstructed from having concordance unknotting number one by Theorem \ref{thm:uc=1}, with the exception of $10_{158}$, which is not obstructed from being concordant to a knot with determinant $5s^2$ and unknotting number one.   
The correction terms of the double branched cover of $10_{158}$ were computed using the method from \cite{osu1}.

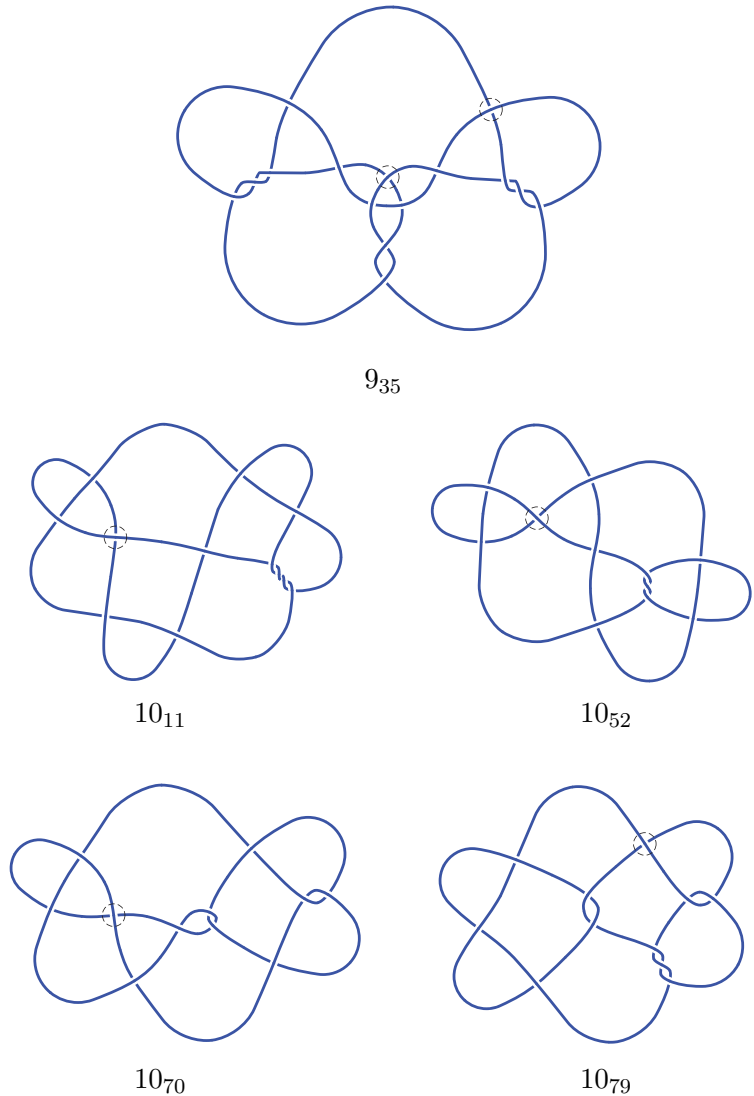
\begin{figure}[htbp]
\begin{center}
\ifpic{
\labellist
\small\hair 3pt
\pinlabel $9_{35}$ at 240 420
\pinlabel $10_{11}$ at 100 210
\pinlabel $10_{52}$ at 380 210
\pinlabel $10_{70}$ at 100 -20
\pinlabel $10_{79}$ at 380 -20
\endlabellist
\includegraphics[scale=0.6]{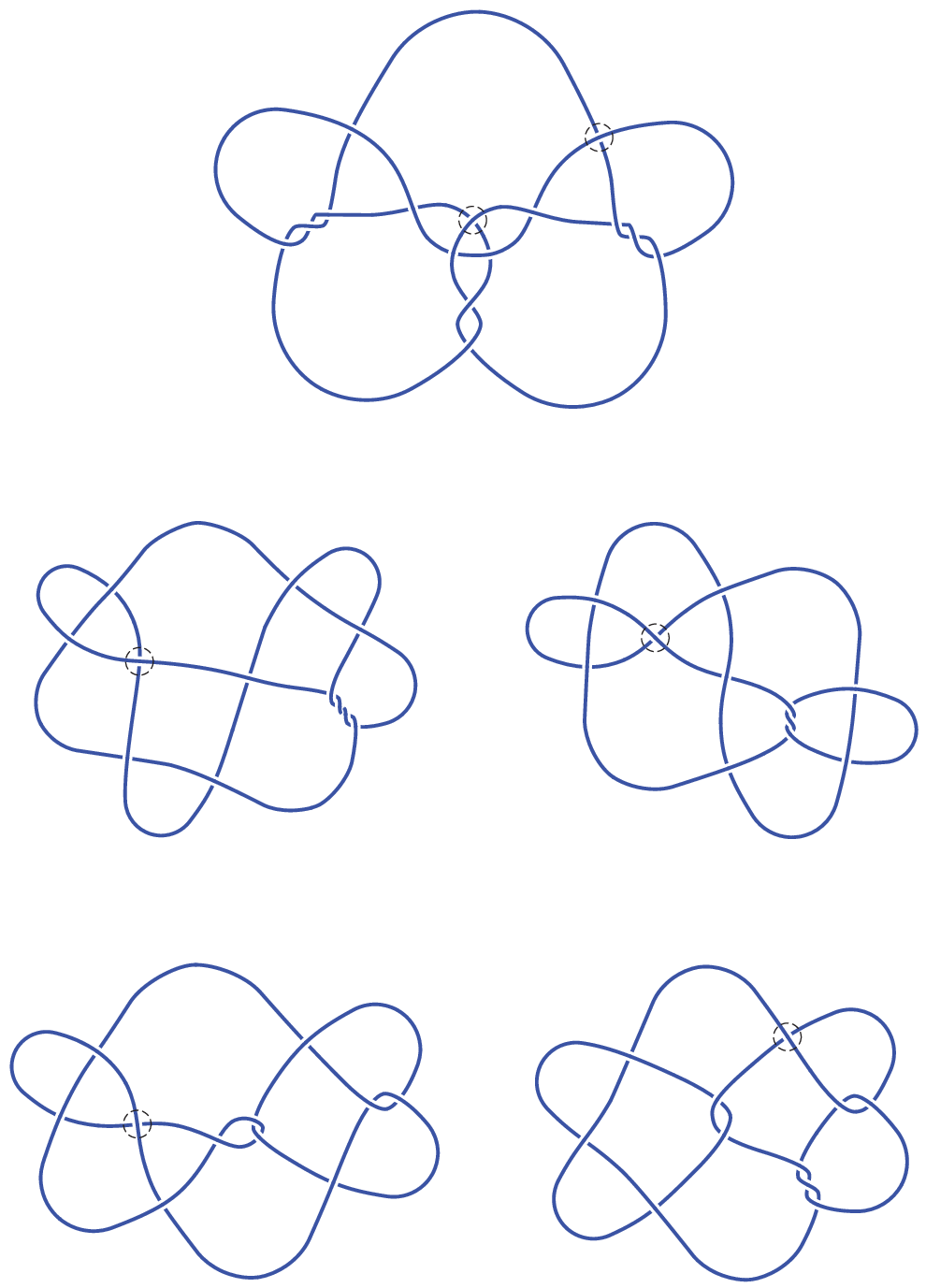}
\vskip 5mm}
\else \vskip 5cm \fi
\begin{narrow}{0.3in}{0.3in}
\caption
{\bf Crossing changes in nonminimal diagrams giving slice knots.}
\label{fig:moreknots}
\end{narrow}
\end{center}
\end{figure}

For all knots in Table \ref{table:10cross} with $u>c^*$ (or $u$ unknown), we have found an appropriate set of crossing changes to convert to a slice knot.  For the most part these are exhibited in the minimal diagram listed in \cite{knotinfo}; the exceptions are shown in Figure \ref{fig:moreknots}.

\begin{table}[!ht]
\tiny
\begin{center}
\begin{tabular}{||c|c|c|c|c|c|c|c||c|c|c|c|c|c|c|c||}
\hline
\rule{0pt}{4ex} 
\rule[-2.3ex]{0pt}{0pt}
Knot  & $\det$ & $u$ & $u_c$ & $u_s$ & $c^*$ & $g^*$ & $\dfrac{|\sigma|}2$ & Knot & $\det$ & $u$ & $u_c$ & $u_s$ & $c^*$ & $g^*$ & $\dfrac{|\sigma|}2$\\
\hline
$10_4$ & $27$ & $2$ & $2$ & $1$ & $1$ & $1$ & $1$ & $10_{79}$ & $61$ & $[2,3]$ & $[2,3]$ & $1$ & $1$ & $1$ & $0$ \\
\hline
$10_6$ & $37$ & $3$ & $[2,3]$ & $2$ & $2$ & $2$ & $2$ & $10_{81}$ & $85$ & $2$ & $2$ & $1$ & $1$ & $1$ & $0$ \\
\hline
$10_{11}$ & $43$ & $[2,3]$ & $[2,3]$ & $1$ & $1$ & $1$ & $1$ & $10_{83}$ & $83$ & $2$ & $2$ & $1$ & $1$ & $1$ & $1$ \\
\hline
$10_{12}$ & $47$ & $2$ & $2$ & $1$ & $1$ & $1$ & $1$ & $10_{86}$ & $85$ & $2$ & $2$ & $2$ & $2$ & $1$ & $0$ \\
\hline
$10_{13}$ & $53$ & $2$ & $2$ & $1$ & $1$ & $1$ & $0$ & $10_{89}$ & $99$ & $2$ & $2$ & $1$ & $1$ & $1$ & $1$ \\
\hline
$10_{15}$ & $43$ & $2$ & $2$ & $1$ & $1$ & $1$ & $1$ & $10_{90}$ & $77$ & $2$ & $2$ & $1$ & $1$ & $1$ & $0$ \\
\hline
$10_{16}$ & $47$ & $2$ & $2$ & $1$ & $1$ & $1$ & $1$ & $10_{93}$ & $67$ & $2$ & $2$ & $1$ & $1$ & $1$ & $1$ \\
\hline
$10_{19}$ & $51$ & $2$ & $2$ & $2$ & $2$ & $1$ & $1$ & $10_{94}$ & $71$ & $2$ & $2$ & $1$ & $1$ & $1$ & $1$ \\
\hline
$10_{20}$ & $35$ & $2$ & $2$ & $2$ & $2$ & $1$ & $1$ & $10_{96}$ & $93$ & $2$ & $2$ & $1$ & $1$ & $1$ & $0$ \\
\hline
$10_{24}$ & $55$ & $2$ & $2$ & $2$ & $2$ & $1$ & $1$ & $10_{97}$ & $87$ & $2$ & $2$ & $2$ & $2$ & $1$ & $1$ \\
\hline
$10_{28}$ & $53$ & $2$ & $2$ & $1$ & $1$ & $1$ & $0$ & $10_{100}$ & $65$ & $[2,3]$ & $[2,3]$ & $2$ & $2$ & $2$ & $2$ \\
\hline
$10_{29}$ & $63$ & $2$ & $2$ & $1$ & $1$ & $1$ & $1$ & $10_{101}$ & $85$ & $3$ & $3$ & $3$ & $3$ & $2$ & $2$ \\
\hline
$10_{34}$ & $37$ & $2$ & $2$ & $1$ & $1$ & $1$ & $0$ & $10_{103}$ & $75$ & $3$ & $1$ & $1$ & $1$ & $1$ & $1$ \\
\hline
$10_{36}$ & $51$ & $2$ & $2$ & $2$ & $2$ & $1$ & $1$ & $10_{105}$ & $91$ & $2$ & $2$ & $2$ & $2$ & $1$ & $1$ \\
\hline
$10_{37}$ & $53$ & $2$ & $2$ & $1$ & $1$ & $1$ & $0$ & $10_{106}$ & $75$ & $2$ & $1$ & $1$ & $1$ & $1$ & $1$ \\
\hline
$10_{38}$ & $59$ & $2$ & $2$ & $1$ & $1$ & $1$ & $1$ & $10_{108}$ & $63$ & $2$ & $2$ & $1$ & $1$ & $1$ & $1$ \\
\hline
$10_{40}$ & $75$ & $2$ & $1$ & $1$ & $1$ & $1$ & $1$ & $10_{109}$ & $85$ & $2$ & $2$ & $2$ & $2$ & $1$ & $0$ \\
\hline
$10_{41}$ & $71$ & $2$ & $2$ & $1$ & $1$ & $1$ & $1$ & $10_{110}$ & $83$ & $2$ & $2$ & $1$ & $1$ & $1$ & $1$ \\
\hline
$10_{43}$ & $73$ & $2$ & $2$ & $1$ & $1$ & $1$ & $0$ & $10_{112}$ & $87$ & $2$ & $2$ & $1$ & $1$ & $1$ & $1$ \\
\hline
$10_{45}$ & $89$ & $2$ & $2$ & $1$ & $1$ & $1$ & $0$ & $10_{115}$ & $109$ & $2$ & $2$ & $1$ & $1$ & $1$ & $0$ \\
\hline
$10_{47}$ & $41$ & $[2,3]$ & $[2,3]$ & $2$ & $2$ & $2$ & $2$ & $10_{116}$ & $95$ & $2$ & $2$ & $2$ & $2$ & $1$ & $1$ \\
\hline
$10_{51}$ & $67$ & $[2,3]$ & $[2,3]$ & $1$ & $1$ & $1$ & $1$ & $10_{117}$ & $103$ & $2$ & $2$ & $1$ & $1$ & $1$ & $1$ \\
\hline
$10_{52}$ & $59$ & $2$ & $2$ & $1$ & $1$ & $1$ & $1$ & $10_{120}$ & $105$ & $3$ & $3$ & $3$ & $3$ & $2$ & $2$ \\
\hline
$10_{53}$ & $73$ & $3$ & $3$ & $3$ & $3$ & $2$ & $2$ & $10_{121}$ & $115$ & $2$ & $2$ & $2$ & $2$ & $1$ & $1$ \\
\hline
$10_{54}$ & $47$ & $[2,3]$ & $[2,3]$ & $1$ & $1$ & $1$ & $1$ & $10_{122}$ & $105$ & $2$ & $2$ & $2$ & $2$ & $1$ & $0$ \\
\hline
$10_{57}$ & $79$ & $2$ & $2$ & $1$ & $1$ & $1$ & $1$ & $10_{125}$ & $11$ & $2$ & $2$ & $1$ & $1$ & $1$ & $1$ \\
\hline
$10_{58}$ & $65$ & $2$ & $2$ & $1$ & $1$ & $1$ & $0$ & $10_{126}$ & $19$ & $2$ & $2$ & $1$ & $1$ & $1$ & $1$ \\
\hline
$10_{61}$ & $33$ & $[2,3]$ & $[2,3]$ & $2$ & $2$ & $2$ & $2$ & $10_{130}$ & $17$ & $2$ & $2$ & $1$ & $1$ & $1$ & $0$ \\
\hline
$10_{64}$ & $51$ & $2$ & $2$ & $1$ & $1$ & $1$ & $1$ & $10_{135}$ & $37$ & $2$ & $2$ & $1$ & $1$ & $1$ & $0$ \\
\hline
$10_{65}$ & $63$ & $2$ & $1$ & $1$ & $1$ & $1$ & $1$ & $10_{138}$ & $35$ & $2$ & $2$ & $1$ & $1$ & $1$ & $1$ \\
\hline 
$10_{67}$ & $63$ & $2$ & $1$ & $1$ & $1$ & $1$ & $1$ & $10_{144}$ & $39$ & $2$ & $2$ & $2$ & $2$ & $1$ & $1$ \\
\hline
$10_{68}$ & $57$ & $2$ & $2$ & $2$ & $2$ & $1$ & $0$ & $10_{148}$ & $31$ & $2$ & $2$ & $1$ & $1$ & $1$ & $1$ \\
\hline
$10_{69}$ & $87$ & $2$ & $2$ & $2$ & $2$ & $1$ & $1$ & $10_{151}$ & $43$ & $2$ & $2$ & $1$ & $1$ & $1$ & $1$ \\
\hline
$10_{70}$ & $67$ & $2$ & $2$ & $1$ & $1$ & $1$ & $1$ & $10_{158}$ & $45$ & $2$ & $[1,2]$ & $1$ & $1$ & $1$ & $0$ \\
\hline
$10_{74}$ & $63$ & $2$ & $1$ & $1$ & $1$ & $1$ & $1$ & $10_{162}$ & $35$ & $2$ & $2$ & $1$ & $1$ & $1$ & $1$ \\
\hline
$10_{76}$ & $57$ & $[2,3]$ & $[2,3]$ & $2$ & $2$ & $2$ & $2$ & $10_{163}$ & $51$ & $2$ & $2$ & $2$ & $2$ & $1$ & $1$ \\
\hline
$10_{77}$ & $63$ & $[2,3]$ & $1$ & $1$ & $1$ & $1$ & $1$ & $10_{165}$ & $39$ & $2$ & $2$ & $2$ & $2$ & $1$ & $1$ \\
\hline
\end{tabular}
\vskip5mm
\begin{narrow}{0.3in}{0.3in}
\caption{
{\bf Invariants of knots with 10 crossings.}  Only knots for which calculation of $u_s$, $u_c$ or $c^*$ is nontrivial are listed.}
\label{table:10cross}
\end{narrow}
\end{center}
\end{table}

\section{Gordian distances between torus knots.}
The Gordian distance $d_G(K_1,K_2)$ between two knots $K_1$ and $K_2$ in the three-sphere is the smallest number of crossing changes required to convert $K_1$ into $K_2$. Note that this notion generalizes the unknotting number, since $u(K)=d_{G}(K,U)$.  The trace of the regular homotopy realizing the crossing changes is a normally immersed annulus in $S^3\times[0,1]$, or in other words an immersed concordance, with one normal double point for each crossing change.  It is natural to consider a weaker measure of distance between knots: we define the \emph{crossing number distance} $d^*(K_1,K_2)$ to be the minimal number of double points in a normally immersed concordance between $K_1$ and $K_2$.  Arguing as in Proposition \ref{prop:movie} this is the same as the \emph{concordance Gordian distance}, in other words the minimum of $d_G(K_1',K_2')$ where $K_i'$ is concordant to $K_i$  for $i=1,2$.  As in the case of embedded concordances, an immersed concordance between $K_1$ and $K_2$ is equivalent to an immersed disk in the four-ball bounded by $-K_1\#K_2$ (one obtains a disk from a concordance by drilling out an arc which avoids the double points, and reverses this by adding a $(3,1)$-handle pair).

There has been a great deal of interest in Gordian distances between torus knots, see for example \cite{borliv,feller}.  The pair $T_{3,10}$ and $T_{5,6}$ seem to be an interesting example since Baader showed in \cite{baader} that they cobound a genus one cobordism in $S^3\times[0,1]$.  As far as we can tell previously known bounds for this pair were $$2\le d_G(T_{3,10},T_{5,6})\le11,$$
with the lower bound coming from Levine-Tristram signatures as detailed below and the upper bound coming from a theorem of Feller \cite[Theorem 2]{feller} which implies that each of these knots has an unknotting sequence of crossing changes with $T_{3,5}$ as an intermediate stage.

\begin{figure}[htbp]
\begin{center}
\ifpic{
\includegraphics[scale=0.8]{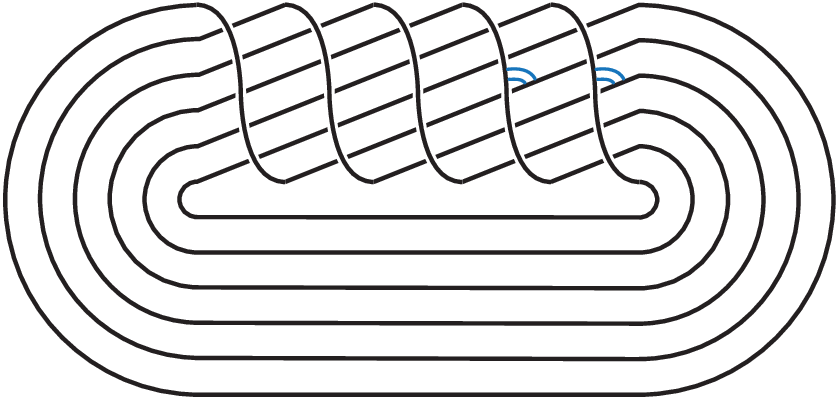}
\vskip 5mm}
\else \vskip 5cm \fi
\begin{narrow}{0.3in}{0.3in}
\caption
{\bf A genus one cobordism between $T_{5,6}$ and $T_{3,10}$.}
\label{fig:torus}
\end{narrow}
\end{center}
\end{figure}

\begin{figure}[htbp]
\begin{center}
\ifpic{
\includegraphics[scale=0.8]{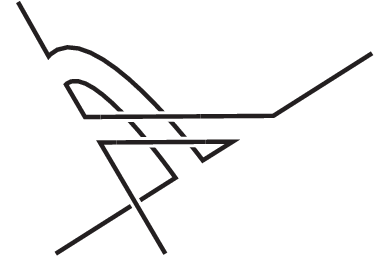}
\vskip 5mm}
\else \vskip 5cm \fi
\begin{narrow}{0.3in}{0.3in}
\caption
{\bf Two band moves realising a crossing change.}
\label{fig:crossingchangebands}
\end{narrow}
\end{center}
\end{figure}

\begin{example}
\label{eg:dGTpq}
Any normally immersed concordance between the torus knots $K_1=T_{3,10}$ and $K_2=T_{5,6}$ has at least three double points, including at least one of each sign.  It is possible to convert $K_1$ into $K_2$ via 5 crossing changes.  Thus the Gordian distance between these knots satisfies
$$3\le d^*(K_1,K_2)\le d_G(K_1,K_2)\le 5.$$
\end{example}
\begin{proof}
For a complex number $z$ of modulus 1 the Levine-Tristram signature $\sigma_z(K)$ and nullity $\eta_z(K)$ of a knot $K$ are defined to be the signature and nullity of $(1-z)V+(1-\overline{z})V^T$, where $V$ is a Seifert matrix for $K$.  The signature $\sigma_z$ agrees for concordant knots $K_1$ and $K_2$ provided $z$ is not a root of the Alexander polynomial of either knot \cite{levine}.  Moreover, the nullity $\eta_z(K)$ vanishes if $z$ is not a root of the Alexander polynomial of $K$.
Also the sum $\sigma_z+\eta_z$ is unchanged or increases (respectively decreases) by two if a positive (resp. negative) crossing is changed (see for example \cite[\S4.2]{borodzikfriedl}).
Using \cite{litherland} we compute Levine-Tristram signatures and nullities of $K_1$ and $K_2$ finding
$$\sigma_{-1}(K_1)=\sigma_\zeta(K_1)=-14,\quad\sigma_{-1}(K_2)=-16,\quad\sigma_\zeta(K_2)=-12,$$
where $\zeta=e^{4\pi i/5}$; moreover, these values are not roots of the Alexander polynomial for either of $K_1$ or $K_2$.   It follows from Proposition \ref{prop:movie} that at least one double point of each sign is required in a normally immersed concordance between $K_1$ and $K_2$.  We will show that Theorem \ref{thm:clasp} obstructs the possibility of such a concordance with two double points.

We let $K=K_1\#-K_2$ and suppose that $K$ bounds a normally immersed disk in the four-ball with two double points.  By the signature data mentioned above, there must be one double point of each sign.  The double branched cover of the torus knot $T_{p,q}$ is the Brieskorn manifold $M(2,p,q)$ \cite{milnor}.  This is the boundary of a negative-definite plumbing tree according to \cite[Theorems 2.1 and 5.1]{neumannraymond}, which is sharp according to results from \cite{plumbed}, enabling us to compute the correction term invariants of these Brieskorn manifolds and hence by additivity of $Y=\Sigma(K)$.

We find that the maximal value taken by the $d$-invariant on $\spinc(Y)$ is $11/10$.  Noting that the determinant of $K$ is 15 which is square-free, there are two possible forms $Q_X$ as in Theorem \ref{thm:clasp} to consider, namely
those represented by
$$Q_1=\left(\begin{matrix}
3 & 0 & 1 & 0\\
0 & 2 & 0 & 1\\
1 & 0 & 2 & 0\\
0 & 1 & 0 & 2
\end{matrix}\right),\quad
Q_2=\left(\begin{matrix}
8 & 0 & 1 & 0\\
0 & 1 & 0 & 1\\
1 & 0 & 2 & 0\\
0 & 1 & 0 & 2
\end{matrix}\right).$$
Neither of these can be the intersection form of a smooth four-manifold whose boundary has 15 \spinc structures, one of whose $d$-invariants is $11/10$; this can be shown using \cite[Theorems 1.2, 9.6]{os4}.  Similar examples are worked out in more detail in \cite{unknotting}.

Finally we observe following Baader \cite{baader} that the band moves in Figure \ref{fig:torus} convert between $K_1$ and $K_2$.  We perform an isotopy on the resulting diagram of $K_2$ by ``sliding a band'': move both ``ends'' of one of the bands past the next band (requiring 4 crossing changes to get past) and then one full revolution around the diagram so that the bands sit on a subdiagram as shown in Figure \ref{fig:crossingchangebands}.  Simplifying this subdiagram and applying one further crossing change gives us back the standard diagram of $K_1$.  Keeping track of signs, we see that one may convert $K_1$ to $K_2$ by changing 2 positive crossings to negative and 3 negative crossings to positive.
\end{proof}
The method of ``sliding bands'' to get upper bounds on Gordian distance may be applied to many of Baader's cobordisms.  Indeed it may be used to show that
\begin{equation}
\label{eqn:dGTpq}
d_G(T_{2a+1,4a+6},T_{2a+3,4a+2})\le 4a+1
\end{equation}
for all $a\ge1$.  To see this, draw a diagram of $T_{2a+3,4a+2}$ as in Figure \ref{fig:torus}, as the closure of the braid on $4a+2$ strands given by a ``$(2a+3)/(4a+2)$ twist'': the top strand passes over all of the other $4a+1$ strands, and this is repeated $2a+3$ times.  Add two band moves to resolve the middle crossing in each of the rightmost two sets of $4a+1$ crossings, again as in Figure \ref{fig:torus};  following Baader, these bands convert the diagram to one of $T_{2a+1,4a+6}$.  Sliding the band on the left $a$ times around the diagram, in a clockwise direction, leads to both bands appearing together as in Figure \ref{fig:crossingchangebands}.
This process involves $4a$ crossing changes: each time the band goes around the torus it has to ``pass through'' the other band, requiring two crossing changes of each sign.
The bands can then be removed using a single crossing change; in total we see that one may obtain $T_{2a+3,4a+2}$ from $T_{2a+1,4a+6}$ by changing a total of $4a+1$ crossings: $2a$ positive to negative, and $2a+1$ negative to positive.
By comparison the results of Feller \cite{feller} give upper bounds which are quadratic in $a$.  For each pair of torus knots listed in \eqref{eqn:dGTpq}
we expect the method of Example \ref{eg:dGTpq} to give the same lower bound of 3 for the crossing number distance, so the range of possibilities is growing linearly in $a$.

\clearpage

\bibliographystyle{amsplain}
\bibliography{concu}

\end{document}